\documentclass[final,3p,times]{elsarticle}

\usepackage{amsmath,amsfonts,amssymb,amscd,subeqnarray,amsthm}
\usepackage{graphicx}
\usepackage{subfig,comment}
\usepackage{siunitx}

\usepackage{stmaryrd}
\usepackage{listings}
\usepackage{boldfonts}
\usepackage{symbols}
\usepackage{multirow}
\usepackage{cases}
\usepackage{tikz, xcolor}
\usepackage{algorithm}
\usepackage{algpseudocode}
\usepackage{natbib}

\usepackage{tikz}
\usetikzlibrary{shapes.geometric}
\usetikzlibrary{shapes,arrows}
\usetikzlibrary{plotmarks}

\tikzstyle{block} = [rectangle, draw, text centered, rounded corners, 
                     text width=7em,    
                     minimum height=3em, node distance=3.3cm]
\tikzstyle{line} = [draw, -latex']
\tikzstyle{cloud} = [draw, ellipse, node distance=2.5cm, minimum height=3em]

\newtheorem{remark}{Remark}

\newcommand{\bdf}[2]{{\textup{\textsf{BDF}}_{#2}\left({#1}\right)}}

\newcommand{\vertiii}[1]{{\left\vert\kern-0.25ex\left\vert\kern-0.25ex\left\vert #1 
    \right\vert\kern-0.25ex\right\vert\kern-0.25ex\right\vert}}

\begin{document}

\begin{frontmatter}
\title{Adaptive enriched Galerkin methods for miscible displacement problems with entropy residual stabilization} 

\author[ut]{Sanghyun Lee\corref{cor1}}
\ead{shlee@ices.utexas.edu}

\author[ut]{Mary F. Wheeler}
\ead{mfw@ices.utexas.edu}

\cortext[cor1]{Corresponding author}
\address[ut]{The Center for Subsurface Modeling, The Institute for Computational Engineering and Sciences,\\ The University of Texas at Austin, Austin, TX 78712}

\begin{abstract}
We present a novel approach to the simulation of miscible displacement by employing adaptive enriched Galerkin finite element methods (EG) coupled with entropy residual stabilization for transport.
In particular, 
numerical simulations of viscous fingering instabilities in heterogeneous porous media and Hele-Shaw cells are illustrated. 
EG is formulated by enriching the conforming continuous Galerkin finite element method (CG) with piecewise constant functions. 
The method provides locally and globally conservative fluxes, which is crucial for coupled flow and transport problems. Moreover, EG has fewer degrees of freedom in comparison with discontinuous Galerkin (DG) and an efficient flow solver has been derived which allows for higher order schemes. Dynamic adaptive mesh refinement is applied in order to save computational cost for large-scale three dimensional applications. In addition, entropy residual based stabilization for high order EG transport systems prevents any spurious oscillations. Numerical tests are presented to show the capabilities of EG applied to flow and transport.
 \end{abstract}


\end{frontmatter}



\section{Introduction}\label{sec:intro}

%
%

Miscible displacement of one 
fluid by another
in a porous medium has attracted considerable attention in subsurface modeling with emphasis on enhanced oil recovery applications~\cite{peaceman1977fundamentals,russell1983finite, wheeler1988numerical}.
Here flow instabilities arising when a fluid with higher mobility 
displaces another fluid with lower mobility is referred to as viscous fingering. 
The latter has been the topic of major physical and mathematical studies for over half a century
\cite{
saffman1958penetration, 
Homsy:1987cm,
Tchelepi:1993cq,
PhysRevLett.109.264503,
malhotra2015experimental,
FLM:376508, 
MIKELIC1991186,
bischofberger2014fingering}.
Recently, viscous fingering has been applied 
for proppant-filled hydraulic fracture propagation~\cite{blyton2015comprehensive,malhotra2014proppant, Lee2016}  to efficiently  transport the proppant to the tip of fractures.

The governing mathematical system that represents the displacement of the fluid mixtures consists of pressure, velocity, and concentration. 
Examples of numerical schemes for approximating this system include the following; 
continuous Galerkin 
\cite{EwingWheeler1980,DOUGLAS1984131,ewing1984galerkin}, 
interior penalty Galerkin \cite{WheelerM-1978aa,wheeler1980interior}, 
finite differences \cite{douglas1982numerical}, 
finite volumes \cite{jenny2006adaptive},
modified method of characteristics \cite{Douglas1982,ewing1983simulation},
mixed finite elements \cite{EWING198473,Douglas1983,darlow1984mixed},
and characteristic-mixed finite elements \cite{russell1986large,ArbotastWheeler1995}.
One of the effective approaches that deals robustly with general partial differential equations  
as well as with equations whose 
type changes within the computational domain 
such as from advection dominated to diffusion dominated  is discontinuous Galerkin (DG)
\cite{LiRiviere2016Fingering,sun2005,Sun2006,Scovazzi201386,Scovazzi2013373}.
DG is well suited for multi-physics applications and for problems with highly varying material properties \cite{riviere2000discontinuous,CNM:CNM464}. 
Combining mixed finite  elements and discontinuous Galerkin was studied  in  \cite{sun2002combined,Li2015107}.

There are three major issues with the above numerical approximations for coupling flow and transport; i) local mass balance, ii) local grid adaptivity, and iii) efficient solution algorithms for  Darcy flow. 
It is well known that
differentiating numerical approximations to obtain a flux suffers from loss of accuracy and the lack of local conservation on  the existing mesh as well as yielding non-physical results for transport with this given flux.
It is important to choose a numerical approximation which preserves local conservation to avoid spurious sources \cite{KAASSCHIETER1995277}.
In addition, the complexities in implementing dynamic grid adaptations can limit the extension of schemes to realistic physical applications. 
Also methods which are computationally costly due to the number of degrees of freedom and lack of efficient solvers prevent developments of higher order methods for large-scale multi-physics problems with highly varying material properties.

In this paper, we introduce a new method for a flow and transport system,
the enriched Galerkin finite element method (EG). This approach provides a locally and globally conservative flux and preserves local mass balance for transport. 
EG is constructed by enriching the
conforming continuous Galerkin finite element method (CG) with piecewise constant functions \cite{BecBurHansLar2003,sunliu2009,LeeLeeWhi15},
with the same bilinear forms as the interior penalty DG schemes.
However, EG has substantially fewer degrees of freedom in comparison with DG 
and a fast effective solver whose cost is roughly that of CG and which can handle an arbitrary order of approximation \cite{LeeLeeWhi15}.
An additional advantage of EG is that 
only those subdomains that require local conservation need to be enriched with a treatment of high order non-matching grids. 

Our high order EG transport system is coupled with an entropy viscosity residual stabilization method introduced in \cite{guermond2011entropy} to 
avoid spurious oscillations near shocks.
Instead of using limiters and non-oscillatory reconstructions,  
this method employs the local residual of an entropy equation to construct the numerical diffusion, which is added as a nonlinear dissipation to the numerical discretization of the system. 
The amount of numerical diffusion added is proportional to the computed entropy residual. 
This technique is independent of mesh and order of approximation and  has been shown to be efficient and stable in solving many physical problems with CG \cite{Guermond2015,FLD:FLD4071,MR3167449, lee_thesis,JL2016}  and DG \cite{Zingan:2013bb}. 

{
In our numerical examples, we illustrate that it is crucial to have dynamic mesh adaptivity in order to reduce computational costs for large-scale three dimensional applications. Earlier work on adaptive local grid refinement in a variety contexts for flow and transport in porous media includes \cite{sun2006dynamic,DOUGLAS1984131,sun2005l2,sun2006posteriori,hornung1997adaptive,edwards1996higher,dahle1992characteristic}.
}
In this paper, we employ the entropy residual for dynamic adaptive mesh refinement to capture the moving interface between the miscible fluids.  
It is shown in \cite{andrews1998posteriori,puppo2004numerical} that the entropy residual can be used as a posteriori error indicator.
Entropy residuals converge to the Dirac measures supported in the shocks as the discretization mesh size goes to zero whereas the residual of the equation converges to zero based on consistency \cite{guermond2011entropy}.
Therefore the entropy residual is able to capture shocks more robustly than general residuals.

In summary, the novelties of the present paper are that 
we establish efficient and robust enriched Galerkin (EG) approximations for miscible displacement problems. 
We couple the high order entropy viscosity stabilization to an EG transport system and implement dynamic mesh adaptivity. In addition, we provide numerical examples to assess the performance of our scheme including viscous fingering instabilities.

The paper is organized as follows. 
The mathematical model is presented
in Section \ref{sec:Model}.  
In Section \ref{sec:discret}, we formulate EG for flow and transport system with the entropy viscosity stabilization method and a global solution algorithm.
Various numerical examples are reported in Section
\ref{sec:num_ex}.

\section{Mathematical Model}\label{sec:Model}

Let $\Omega\subset \Reals{d}$ be a bounded polygon (for $d=2$) or polyhedron (for $d=3$) with Lipschitz boundary $\partial \Omega$ and $(0,\mathbb{T}]$ is the computational time interval with $\mathbb{T} > 0$.
We consider a multi-component miscible displacement system with a single phase slightly compressible flow.
The advection-diffusion transport system for the miscible components $i$ is given as 
\begin{equation}
\dfrac{\partial}{\partial t}( \varphi \rho  c_i) + \nabla \cdot (\rho \bu c_i - \varphi \rho \bD(\bu) \nabla c_i)  = \tilde{q}_i , \mbox{ in } \Omega \times (0,\mathbb{T}],   \label{eqn:main_transport_i}
\end{equation}
where 
$\varphi$ is the porosity, 
 $\bu : \Omega \times [0,\mathbb{T}] \rightarrow \mathbb{R}^d$ is the velocity, 
$c_i : \Omega \times (0,T] \rightarrow \mathbb{R}$  is the advected mass fraction of the component $i$ of the solution,
and the average density $\rho$ is defined as 
\begin{equation}
\rho := \left(\sum_{i=1}^{N_c} \dfrac{c_i}{\rho_i} \right)^{-1}
\end{equation}
with the total number of components $N_c$ 
by assuming there is no volume change in mixing.
For convenience, we assume only two components in our case ($N_c=2$ and $i=1,2$), in particular we set $c:=c_1$ and $1-c := c_2$.
This leads to solving for only one component as follows;
\begin{equation}
\dfrac{\partial}{\partial t}( \varphi \rho c) + \nabla \cdot (\rho \bu c - \varphi \rho \bD(\bu) \nabla c)  = \tilde{q} , \mbox{ in } \Omega \times (0,\mathbb{T}],   
\label{eqn:main_transport_0}
\end{equation}
where $\tilde{q} := \tilde{q}_1$ without loss of generality 
and the remaining component is obtained by the relation $c_1 + c_2 = 1$.
Since the flow is assumed to be slightly compressible, the compressibility coefficient satisfies $c_F^i \ll 1 $ in the relationship 
$$
\rho_i(p) \approx \rho^i_0(1 + c_F^i p),
$$
where $\rho^i_0$ is the initial density of the fluid for each component. 
For simplicity, we assume that
each component has the same initial density ($\rho_0 := \rho_0^1 = \rho_0^2$) and compressibility coefficient ($c_F := c_F^1 = c_F^2$). 
Under above assumptions, we can rewrite \eqref{eqn:main_transport_0} to 
\begin{equation}
\dfrac{\partial}{\partial t}( \varphi \rho_0  c) + \nabla \cdot (\rho_0 \bu c - \varphi \rho_0 \bD(\bu) \nabla c)  = \tilde{q} , \mbox{ in } \Omega \times (0,\mathbb{T}].
\label{eqn:main_transport}
\end{equation}
Here  $\tilde{q} := \tilde{c} q$ and $\tilde{c}, q$ are the concentration source/sink term and flow source/sink term, respectively.  
If $q>0$, $\tilde{c}$ is the injected concentration $c_q$ and 
if $q<0$, $\tilde{c}$ is the resident concentration $c$.
The dispersion/diffusion tensor is defined as,
\begin{equation}
\bD(\bu) := d_m \bI + | \bu | \left( \alpha_l \bE(\bu) + \alpha_t ( \bI - \bE(\bu)) \right),
\label{eqn:diff_dispersion}
\end{equation}
where 
$$
(\bE(\bu))_{ij} := 
\dfrac{(\bu_i \bu_j)}{| \bu |^2}, 1\leq i,j \leq d
$$ 
is the tensor that projects onto the $\bu$ direction, 
$d_m >0$ is the molecular diffusivity,
$\alpha_l >0$ is the longitudinal, 
and $\alpha_t >0$ is the transverse dispersivities \cite{wheeler1988numerical}.

Next, the flow is described by following
\begin{equation}
\dfrac{\partial}{\partial t}(\varphi \rho_0 c_F p  ) + \nabla\cdot \left ( \rho_0 \bu \right ) = q  \mbox{ in } \Omega \times (0,\mathbb{T}], \label{eqn:main_pressure} 
\end{equation}
where $p:\Omega \times [0,\mathbb{T}] \rightarrow \mathbb{R}$ is the fluid pressure.
The velocity $\bu : \Omega \times [0,\mathbb{T}] \rightarrow \mathbb{R}^d$ is defined by Darcy's law 
\begin{equation}
\bu = -\dfrac{K}{\mu(c)}  \big( \nabla p -\rho \bg \big) ,  \mbox{ in } \Omega \times (0,\mathbb{T}],
\label{eqn:main_velocity}
\end{equation}    
where $K := K(\bx) = (K^{ij}(\bx))_{i,j=1,\cdots,d}$ for $\bx \in \Omega$,  denotes the permeability coefficient in $[L^\infty(\Omega)]^{d\times d}$ and 
the function $\mu := \mu(c)$ is the fluid viscosity in $L^\infty(\Omega)$, both of which may have jump discontinuities. 
We assume that $K$ is uniformly symmetric positive definite, with respect to an initial non-overlapping (open) subdomain partition of the domain $\Omega$.
Set $\mathcal{T}_S = \{\Omega_m\}_{m=1}^{M},$ with $\cup_{m=1}^M \overline{\Omega}_m=\overline{\Omega}$ and ${\Omega}_m \cap {\Omega}_n =\emptyset$ for $n\neq m$.  
The (polygonal or polyhedral) regions $\Omega _m\;, m=1, \ldots, M,$ may involve complicated geometry. 
We define $\bkappa := \bkappa(c) := K/\mu(c) $ and
the gravity field is denoted by $\bg$ but neglected in our discussion.
%

\subsection{Boundary and Initial conditions}

The boundary of $\Omega$ for transport system, denoted by $\partial \Omega$,
 is decomposed into two parts 
$\Gamma_{\rm in}$ and $\Gamma_{\rm out}$, the inflow and outflow boundary, respectively, (i.e. $\overline{\partial \Omega} = \overline{\Gamma}_{\rm in}  \cup \overline{\Gamma}_{\rm out} $.
) Those are defined as 
\begin{equation}
\Gamma_{\rm in} := \{ \bx \in \partial \Omega : \bu\cdot\bn < 0\} \quad \mbox{ and } \quad \Gamma_{\rm out} := \{ \bx \in \partial \Omega : \bu\cdot\bn \geq 0\},
\label{eqn:main_bd_def}
\end{equation}
where $\bn$ denotes the unit outward normal vector to $\partial \Omega$. 
For each boundary, we employ the following boundary conditions
\begin{align}
(\rho_0 \bu c - \varphi \rho_0 \bD(\bu) \nabla c)  \cdot \bn = 
c_\text{in} \rho_0 \bu \cdot \bn, 
& \mbox{ on } \Gamma_{\text{in}} \times (0,\mathbb{T}], \label{main_bd_c_in} \\
(-\varphi \rho_0  \bD(\bu)\nabla c) \cdot \bn = 0 , 
& \mbox{ on } \Gamma_{\text{out}} \times (0,\mathbb{T}],
\label{main_bd_c_out}
\end{align}
where $c_{\text{in}}$ is a given inflow boundary value.

For the flow problem, the boundary is decomposed into two parts $\Gamma_D$ and $\Gamma_N$ so that $\overline{\partial \Omega} = \overline{\Gamma}_D \cup \overline{\Gamma}_N$ and we impose
\begin{align}
p = g_{_D} & \mbox{ on } \Gamma_D \times (0,\mathbb{T}], \quad \\ 
\rho_0 \bu \cdot \bn = g_{_N} & \mbox{ on } \Gamma_N \times (0,\mathbb{T}],  
\end{align}
where $g_D \in L^2(\Gamma_D)$ and $g_N \in L^2(\Gamma_N)$ are the each Dirichlet and Neumann boundary conditions, respectively.

The above systems are supplemented by initial conditions 
$$
c(\bx,0) = c_0(\bx), \mbox{ and }
p(\bx,0) = p_0(\bx), \quad \forall \bx \in \Omega.
$$
\section{Numerical Method}
\label{sec:discret}
Let $\mathcal{T}_h$ be the shape-regular (in the sense of Ciarlet)  triangulation by a family of partitions of $\O$ into $d$-simplices $\K$ (triangles/squares in $d=2$ or tetrahedra/cubes in $d=3$). We denote by $h_{\K}$ the diameter of $\K$ and we set $h=\max_{\K \in \Th} h_{\K}$.  
Also we denote by $\Eh$ the set of all edges and by $\Eho$ and $\Ehb$ the collection of all interior and boundary edges, respectively. 
In the following notation, we assume edges for two dimension but the results hold analogously for faces in three dimensional case.
For the flow problem, the boundary edges $\Ehb$ can be further decomposed into $\Ehb = \mathcal{E}_h^{D,\partial} \cup \mathcal{E}_{h}^{N,\partial}$, where $\mathcal{E}_h^{D,\partial}$ is the collection of edges where the Dirichlet boundary condition is imposed, while $\mathcal{E}_h^{N,\partial}$ is the collection of edges where the Neumann boundary condition is imposed. 
In addition, we let  $\mathcal{E}_h^{1} := \Eho \cup \mathcal{E}_h^{D,\partial}$ and $\mathcal{E}_h^{2} := \Eho \cup \mathcal{E}_h^{N,\partial}$. 
For the transport problem,  
the boundary edges $\Ehb$ decompose into $\Ehb = \mathcal{E}_h^{\text{in}} \cup \mathcal{E}_{h}^{\text{out}}$, where $\mathcal{E}_h^{\text{in}}$ is the collection of edges where the inflow boundary condition is imposed, while $\mathcal{E}_h^{\text{out}}$ is the collection of edges where the outflow boundary condition is imposed.

The space $H^{s}(\Th)$ $(s\in \Reals{})$ is the set of element-wise $H^{s}$ functions on $\mathcal{T}_h$, and $L^{2}(\Eh)$ refers to the set of functions whose traces on the elements of $\Eh$ are square integrable. Let $\mathbb{Q}_k(\K)$ denote the space of polynomials of partial degree at most $k$. 
Regarding the time discretization, given an integer $N \geq 2$, we define a partition of the time interval 
$0 =: t^0 < t^1 < \cdots < t^N:= \mathbb{T}$ and denote $\delta t := t^n - t^{n-1}$ for the uniform time step.
Throughout the paper, we use the standard notation for Sobolev spaces \cite{AdamsR-1975aa} and their norms. For example, let $E \subseteq \Omega$, then $\|\cdot\|_{1,E}$ and $|\cdot|_{1,E}$ denote the $H^1(E)$ norm and seminorm, respectively. 
For simplicity, we eliminate the subscripts on the norms if $E = \Omega$.
For any vector space $\bX$, $\bX^d$ will denote the vector space of size d, whose components belong to $\bX$ and $\bX^{d\times d}$ will denote the $d \times d$ matrix whose components belong to $\bX$.

We  introduce the space of piecewise discontinuous polynomials of degree $k$ as
\begin{equation}
M^k(\mathcal{T}_h) := \left \{ \psi \in L^2(\Omega) | \ \psi_{|_{\K}} \in \mathbb{Q}_k(\K), \ \forall \K \in \mathcal{T}_h \right \}, 
\end{equation}
and let $M_0^k(\mathcal{T}_h)$ be the subspace of $M^k(\mathcal{T}_h)$ consisting of continuous piecewise polynomials;
\begin{equation*}
M_0^k(\mathcal{T}_h) = M^k(\mathcal{T}_h) \cap \mathbb{C}_0(\Omega). 
\end{equation*}
The enriched Galerkin finite element space, denoted by $V_{h,k}^{\textsf{EG}}$ is defined as
\begin{equation}
V_{h,k}^{\textsf{EG}}(\mathcal{T}_h)  := M^k_0(\mathcal{T}_h) + M^0(\mathcal{T}_h),
\end{equation}
where $k \geq 1$, also see \cite{BecBurHansLar2003,sunliu2009,LeeLeeWhi15,LeeLeeWhi16a} for more details.

\begin{figure}[!h]
\centering
\subfloat[$CG_1$]{
{
\begin{tikzpicture}
\draw (0,0) -- (2,0); 
\draw (0,1) -- (2,1); 
\draw (2,0) -- (2,2); 
\draw (2,2) -- (0,2); 
\draw (0,2) -- (0,0);
\draw (1,2) -- (1,0); 
\draw[solid,fill=red] (1,1) circle (0.05);
\draw[solid,fill=red] (1,0) circle (0.05);
\draw[solid,fill=red] (0,1) circle (0.05);
\draw[solid,fill=red] (2,1) circle (0.05);
\draw[solid,fill=red] (1,2) circle (0.05);
\draw[solid,fill=red] (0,0) circle (0.05);
\draw[solid,fill=red] (2,0) circle (0.05);
\draw[solid,fill=red] (2,2) circle (0.05);
\draw[solid,fill=red] (0,2) circle (0.05);
\end{tikzpicture}
}
}
\hspace{1cm}
\subfloat[$DG_1$]{
{
\begin{tikzpicture}
\draw (0,0) -- (2,0); 
\draw[dashed] (0,1) -- (2,1); 
\draw (2,0) -- (2,2); 
\draw (2,2) -- (0,2); 
\draw (0,2) -- (0,0);
\draw[dashed] (1,2) -- (1,0); 
\draw[solid,fill=red] (0.9,1.1) circle (0.05);
\draw[solid,fill=red] (1.1,1.1) circle (0.05);
\draw[solid,fill=red] (0.9,0.9) circle (0.05);
\draw[solid,fill=red] (1.1,0.9) circle (0.05);
\draw[solid,fill=red] (0.9,0) circle (0.05);
\draw[solid,fill=red] (1.1,0) circle (0.05);
\draw[solid,fill=red] (0,1.1) circle (0.05);
\draw[solid,fill=red] (0,0.9) circle (0.05);
\draw[solid,fill=red] (2,1.1) circle (0.05);
\draw[solid,fill=red] (2,0.9) circle (0.05);
\draw[solid,fill=red] (1.1,2) circle (0.05);
\draw[solid,fill=red] (0.9,2) circle (0.05);
\draw[solid,fill=red] (0,0) circle (0.05);
\draw[solid,fill=red] (2,0) circle (0.05);
\draw[solid,fill=red] (2,2) circle (0.05);
\draw[solid,fill=red] (0,2) circle (0.05);
\end{tikzpicture}
}
}
\hspace{1cm}
\subfloat[$EG_1$]{
{
\begin{tikzpicture}
\draw (0,0) -- (2,0); 
\draw (0,1) -- (2,1); 
\draw (2,0) -- (2,2); 
\draw (2,2) -- (0,2); 
\draw (0,2) -- (0,0);
\draw (1,2) -- (1,0); 
\draw[solid,fill=red] (1,1) circle (0.05);
\draw[solid,fill=red] (1,0) circle (0.05);
\draw[solid,fill=red] (0,1) circle (0.05);
\draw[solid,fill=red] (2,1) circle (0.05);
\draw[solid,fill=red] (1,2) circle (0.05);
\draw[solid,fill=red] (0,0) circle (0.05);
\draw[solid,fill=red] (2,0) circle (0.05);
\draw[solid,fill=red] (2,2) circle (0.05);
\draw[solid,fill=red] (0,2) circle (0.05);

\node[mark size=2pt,color=blue] at (0.4,0.5) {\pgfuseplotmark{triangle*}};
\node[mark size=2pt,color=blue] at (1.4,0.5) {\pgfuseplotmark{triangle*}};
\node[mark size=2pt,color=blue] at (0.4,1.5) {\pgfuseplotmark{triangle*}};
\node[mark size=2pt,color=blue] at (1.4,1.5) {\pgfuseplotmark{triangle*}};

\end{tikzpicture}
}
}
\caption{
{
Comparison of the degrees of freedom for a two-dimensional Cartesian grid ($\mathbb{Q}_1$) with linear CG, DG and EG approximations. 
Here the triangle in the middle of the grid at (c) indicates a piece-wise constant ($M^0(\mathcal{T}_h)$).
}
}
\label{fig:dof}
\end{figure}
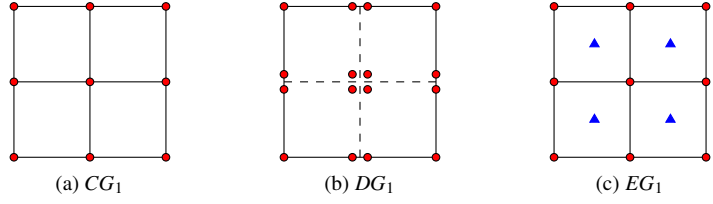

\begin{remark}
We remark that the degrees of freedom for $V_{h,1}^{\textsf{EG}}(\mathcal{T}_h) $ with large enough number of grids 
is approximately one half and {one fourth} the degrees of freedom of the linear DG space, in two and three space dimensions, respectively. {\color{blue}See Figure \ref{fig:dof}}.
\end{remark}

For the EG formulation, we employ the weighted Interior Penalty (IP) methods as presented in \cite{ErnA_StephansenA_ZuninoP-2009aa,StenbergR-1998aa}. First, we define the coefficient $\bkappa_T$ by 
\begin{equation}
\bkappa_T := \bkappa |_T, \quad \forall T \in \mathcal{T}_h.
\end{equation}
Following \cite{ArnoldD_BrezziF_CockburnB_MariniL-2001aa}, for any $e \in \Eho$, let $\K^{+}$ and $\K^{-}$ be two neighboring elements such that  $e = \partial \K^{+}\cap \partial \K^{-}$.
We denote by $h_{e}$ the length of the edges $e$.
Let $\n^{+}$ and $\n^{-}$ be the outward normal unit vectors to  $\partial T^+$ and $\partial T^-$, respectively ($\n^{\pm} :=\n_{|\K^{\pm}}$). 
For any given function $\xi$ and vector function $\bxi$, defined on the triangulation $\mathcal{T}_h$, we denote $\xi^{\pm}$ and $\bxi^{\pm}$ by the restrictions of $\xi$ and $\bxi$ to $T^\pm$, respectively. 
Given certain weight $\delta_e \in [0,1]$, we define the weighted average $\av{\cdot}_{\delta_{e}}$ as follows: for $\zeta \in L^2(\mathcal{T}_h)$ and $\taub \in L^2(\mathcal{T}_h)^d$,
\begin{equation}\label{av-w}
\av{\zeta}_{\delta_{e}} := \delta_e \zeta^+ + (1-\delta_e) \zeta^- 
\quad \mbox{ and } \quad 
\av{\taub}_{\delta_{e}} := \delta_e\taub^+ +   (1- \delta_e) \taub^- \quad \mbox{on } e\in
\Eho.
\end{equation}
The usual average $\av{\cdot}_{1/2}$ will be simply denoted by $\av{\cdot}$,
\begin{align*}
\av{\zeta} := \av{\zeta}_{1/2} \quad \mbox{ and } \quad \av{\taub} := \av{\taub}_{1/2}, \quad \mbox{on } e \in \Eho.
\end{align*}
On the other hand, for $e \in \Ehb$, we set $\av{\zeta}_{\delta_{e}} := \av{\zeta} = \zeta$ and $\av{\taub}_{\delta_{e}} := \av{\taub} = \taub$. 
The jump across the interior edge will be defined as usual: 
\begin{align*}
\jump{\zeta} = \zeta^+\n^++\zeta^-\n^- \quad \mbox{ and } \quad \jtau = \taub^+\cdot\n^+ + \taub^-\cdot\n^- \quad \mbox{on } e\in \Eho. 
\end{align*}
For $e\in \Ehb$, we set $\jump{\zeta} = \zeta\n$. 
The choice of the weights has been investigated in  \cite{Li:2015fl,BurmanE_ZuninoP-2006aa,Bastian:2014fn} and references cited therein. 
In this paper, we consider the following choice of weights in terms of $\bkappa$ given as follows. 
We first define, 
$$
\kappa^+_e := (\bn^+)^{\top}  \bkappa^+  \bn^+
\quad \mbox{ and } \quad 
\kappa^-_e := (\bn^+)^{\top} \bkappa^-  \bn^+ 
$$
with fixed unit normal direction. Then, the weight is  chosen as 
$$
\delta_e = \beta_e := \dfrac{\kappa^-_e}{\kappa^+_e + \kappa^-_e}.
$$
The coefficient $\kappa_e$ is defined as the harmonic mean of $\kappa^{+}_e$ and $\kappa^{-}_e$ by
\begin{equation}\label{defKE}
\kappa_e := \frac{2 \kappa^{+}_e \kappa^{-}_e }{ \kappa^+_e + \kappa^{-}_e}.
\end{equation} 
We note that the weights $\{\delta_e\}_{e\in \Eho}$ depend on the coefficient $\bkappa$ and they may vary over all interior edges. 
We also note that for each $e \in \Eho$ the weighted average $ \av{ \bkappa \nabla  v}_{\beta_e}$ for $\forall v \in H^1(\mathcal{T}_h)$, can be rewritten as 
$$
\av{\bkappa \nabla v}_{\beta_e} 
= 
\beta_e (\bkappa^{+}(\nabla v)^{+})+(1-\beta_{e})
(\bkappa^{-}(\nabla v)^{-}).
$$
For inner products, we use the notations: 
\begin{align*}
&(v,w)_{\Th}:=\dyle\sum_{\K \in \Th} \int_{\K} v\, w dx, \quad \forall\,\, v ,w \in L^{2} (\mathcal{T}_h), \\
&\langle v, w\rangle_{\Eh}:=\dyle\sum_{e\in \Eh} \int_{e} v\, w \,d\gamma, \quad \forall\, v, w \in L^{2}(\Eh).
\end{align*}
{
For example,  $\psi_{\textsf{EG}} \in V_{h,k}^{\textsf{EG}}(\mathcal{T}_h)$ decomposes into 
$\psi_{\textsf{EG}} = \psi_{\textsf{CG}}+  \psi_{\textsf{DG}}$,
where 
$\psi_{\textsf{CG}} \in M^k_0(\mathcal{T}_h)$  
and
$\psi_{\textsf{DG}} \in M^0(\mathcal{T}_h)$. 
Thus the inner product $(\psi_{\textsf{EG}} ,\psi_{\textsf{EG}} ) = 
(\psi_{\textsf{CG}} ,\psi_{\textsf{CG}} ) + (\psi_{\textsf{CG}} ,\psi_{\textsf{DG}} ) + (\psi_{\textsf{DG}} ,\psi_{\textsf{CG}} ) + (\psi_{\textsf{DG}} ,\psi_{\textsf{DG}} )$ creates a block matrix which has the following form,
$$
\begin{pmatrix} 
\psi_{\textsf{CG}} \psi_{\textsf{CG}} & \psi_{\textsf{CG}} \psi_{\textsf{DG}} \\
\psi_{\textsf{DG}} \psi_{\textsf{CG}} & \psi_{\textsf{DG}} \psi_{\textsf{DG}}
\end{pmatrix}.
$$
}
Finally, we introduce the interpolation operator $\Pi_h$ for the space $V_{h,k}^{\textsf{EG}}$ as
\begin{equation}\label{interp}
\Pi_h v = \Pi_0^k v + Q^0 ( v - \Pi_0^k v), 
\end{equation}
where $\Pi_0^k$ is a {continuous} interpolation operator onto the space $M_0^k(\mathcal{T}_h)$, {and} $Q^0$ is the $L^2$ projection onto the space $M^0(\mathcal{T}_h)$.  See \cite{LeeLeeWhi15} for more details.
\subsection{Numerical Approximation of the Pressure System}
The locally conservative EG is used for the space approximation of the pressure system  \eqref{eqn:main_pressure}-\eqref{eqn:main_velocity}. 
We refer to \cite{LeeLeeWhi15} for 
a mathematical discussion on its stability and error convergence properties with an efficient solver which we employ here.

\subsubsection{Space and Time Discretization}
The time discretization is carried out by choosing 
$N\in \mathbb{N}$, the number of time steps, and setting the time step to be $\delta t = \mathbb{T}/N$. We set $t^n = n \delta t$ and for a time dependent function we denote $\phi^n = \phi(t^n)$ and 
$\phi^{\delta t} = \{ \phi^n \}_{n=0}^{N}$.
Over these sequences we define the time stepping  operator
\begin{equation}
  \bdf{\phi^{n+1}}{m} :=  
\begin{cases}
\dfrac{1}{\delta t}    (\phi^{n+1} - \phi^n) & \text{if } m=1, \vspace*{0.06in} \\
\dfrac{1}{2\delta t} \left( 3 \phi^{n+1} - 4\phi^n + \phi^{n-1} \right)   & \text{if } m=2,
\end{cases}
\label{eqn:bdf}
\end{equation}
for different discretization order $m$. First order Euler and backward differentiation formula (BDF2) is applied for $m=1,2$, respectively.
  
The EG finite element space approximation of the pressure $p(\bx,t)$ is denoted by $P(\bx,t) \in V^{\textsf{EG}}_{h,k}(\mathcal{T}_h) $ and 
we let $P^n := P(\bx,t^n)$ for time discretization, $0 \leq n \leq N$.
Here $m=2$ is chosen for order of the time-stepping scheme. 
We set an initial condition for the pressure as  $P^0 := \Pi_h  p(\cdot,0)$. Let  $g_D^{n+1}, g_N^{n+1}$ and $q^{n+1}$ are approximations of $g_D(\cdot, t^{n+1}),g_N(\cdot, t^{n+1})$ and $q(\cdot, t^{n+1})$ on $\Gamma_D$, $\Gamma_N$ and $\Omega$, respectively at time $t^{n+1}$.
Assuming for the time being that $c(\cdot,t^{n+1})$ and 
$\bkappa(t^{n+1}) := \bkappa(c(\cdot,t^{n+1}))$ are known, 
the time stepping algorithm reads as follows: Given $P^{n}$, find
\begin{equation}\label{eq:egscheme}
P^{n+1} \in V_{h,k}^{\textsf{EG}}(\mathcal{T}_h) \mbox{  such that  } \calS_\theta(P^{n+1},w) = \mathcal{F}_\theta(w), \quad \forall\, w\in V_{h,k}^{\textsf{EG}}(\mathcal{T}_h) , \,
\end{equation}
where $\calS_\theta$ and $\mathcal{F}_\theta$ are the bilinear form and linear functional defined as 
\begin{multline}
\calS_\theta(P^{n+1},w) := 
\left( \rho_0  \varphi c_F  \bdf{P^{n+1}}{2}, w \right)_{\mathcal{T}_h}
+ 
\left (\rho_0 \bkappa(t^{n+1}) \nabla P^{n+1},\nabla w \right )_{\Th} - \left \langle \rho_0  \av{ \bkappa(t^{n+1}) \nabla P^{n+1} }_{\beta_{e}}, \jump{w} \right \rangle_{\mathcal{E}_h^{1}} \\
 \quad 
+ 
\theta \left \langle \jump{P^{n+1}}, \rho_0 \av{ \bkappa(t^{n+1}) \nabla w}_{\beta_{e}} \right \rangle_{\mathcal{E}_h^{1}} 
+ \dfrac{\alpha(k)}{h_{e}} \rho_0 \left\langle   \bkappa_e(t^{n+1}) \jump{P^{n+1}},\jump{w} \right\rangle_{\mathcal{E}_h^{1}}, \nonumber 
\end{multline}
and 
\begin{equation}
\mathcal{F}_\theta(w) := 
\left(q^{n+1},w\right)_{\mathcal{T}_h} 
- \left \langle  g_N^{n+1},\jump{w} \right \rangle_{\mathcal{E}_h^{N,\partial}}
+ \theta \left \langle g_D^{n+1}, \rho_0 \av{ \bkappa(t^{n+1}) \nabla w}_{\beta_{e}} \cdot \bn \right \rangle_{\mathcal{E}_h^{D,\partial}} 
+\dfrac{\alpha(k)}{h_{e}}  \left\langle  \rho_0   \bkappa_e(t^{n+1}) g_D^{n+1},\jump{w}\right\rangle_{\mathcal{E}_h^{D,\partial}}.
\label{ftheta}
\end{equation}

Here $\alpha(k)$ is a penalty parameter that can vary {on edges where $k$ is the degree of polynomial employed for the space $V_{h,k}^{\textsf{EG}}$.
The {choice} of $\theta$ leads to different EG algorithms. For example, i) $\theta = -1$ for SIPG($\beta$)$-k$ methods, 
ii) $\theta = 1$ for NIPG($\beta$)$-k$ methods, and iii) $\theta = 0$ for IIPG($\beta$)$-k$ method. 
For this paper, we set $\theta = 0$ for the flow problem to satisfy the compatibility condition 
that implies
if the concentration $(c)$ is identically equal to a positive constant ($c^*$) then $(c = c^*)$ is preserved in transport (see equation  (36) in \cite{DawsonC_SunS_WheelerM-2004aa} for more details). 

\subsubsection{Locally conservative flux}
Such conservative flux variables $\bU^{n+1}$ can be obtained as \cite{sunliu2009} and 
details for conservation analyses and EG approximation estimate for the flux in our problem is discussed in \cite{LeeLeeWhi15}.
Let $P^{n+1}$ be the solution to the \eqref{eq:egscheme}, then  we define the globally and locally conservative flux variables $\bU^{n+1}$ at time step $t^{n+1}$ by the following :
\begin{subeqnarray}\label{flux}
\bU^{n+1} |_{T} &=& - \bkappa(t^{n+1}) \nabla P^{n+1}, \quad  \forall T \in \mathcal{T}_h \\ 
\bU^{n+1} \cdot \bn|_{e} &=& - 
\av{\bkappa(t^{n+1}) \nabla  P^{n+1}} \cdot \bn + \alpha(k) h_e^{-1} \bkappa_e(t^{n+1}) \jump{P^{n+1}}, 
\quad \forall e \in \mathcal{E}_h^I, \\ 
\bU^{n+1} \cdot \bn|_{e} &=& g_N^{n+1}, \quad \forall e \in \mathcal{E}_h^{N,\partial}, \\
\bU^n \cdot \bn|_{e} &=& - \bkappa(t^{n+1}) \nabla P^{n+1} \cdot \bn + \alpha(k) h_e^{-1} \bkappa(t^{n+1}) \left( P^{n+1} - g_D^{n+1} \right ), \quad \forall e \in \mathcal{E}_h^{D,\partial},
\end{subeqnarray}
where $\bn$ is the unit normal vector of the boundary edge $e$ of $T$.

\subsection{Numerical Approximation of the Transport System}
The bilinear form of EG coupled with an entropy residual stabilization is employed for modeling the transport system \eqref{eqn:main_transport} with  higher order approximations.
The time stepping is done by using a second order backward Euler ($m=2$ for \eqref{eqn:bdf}).  
Stability and error convergence analyses for the approximation are provided in \cite{LeeLeeWhi16a}. 

\subsubsection{Space and Time Discretization}
Let $C(\bx,t)$ be the space approximation of the concentration function $c(\bx,t)$ and the time approximation of $C(\bx, t^n), 0 \leq n \leq N$ be denoted by $C^n$.  
We set an initial condition for the concentration as  $C^0 := \Pi_h  C(\cdot,0)$.
The discretized system to find $C^{n+1} \in V_{h,\tilde{k}}^{\textsf{EG}}(\mathcal{T}_h)$ is given as follow: Given $C^{n}$, find
\begin{equation}
C^{n+1} \in V_{h,\tilde{k}}^{\textsf{EG}}(\mathcal{T}_h) \mbox{  such that  }
{\mathcal{A}}(C^{n+1}, v) = {\mathcal{L}}(v), \quad \forall v \in V_{h,\tilde{k}}^{\textsf{EG}}(\mathcal{T}_h), 
\label{eqn:discrete_transport} 
\end{equation}
where,
\begin{multline}
{\cal{A}}(C^{n+1}, v) := 
\left( \varphi \rho_0   \bdf{C^{n+1}}{2}, v \right)_{\mathcal{T}_h}
+
\left( \varphi \rho_0  \bD(\bU^n) \nabla C^{n+1} 
- \rho_0 \bU^{n+1} C^{n+1} ,  \nabla v \right)_{\mathcal{T}_h} \\
-
\left\langle  \varphi  \rho_0  \av{\bD(\bU^n) \nabla C^{n+1} }, \jump{v} \right\rangle_\Eho
+
\langle { (C^{n+1})}^* \rho_0  \bU^{n+1} , \jump{v} \rangle_\Eho
- (c  (q^{n+1})^-, v )_{\mathcal{T}_h}  \\
+
\dfrac{\alpha_{c}(\tilde{k})}{ h_e} \rho_0 \left\langle  \jump{C^{n+1}},\jump{v}\right\rangle_\Eho
+ 
\left\langle C^{n+1}  \rho_0  \bU^{n+1}  \cdot\bn,  v 
\right\rangle_{\mathcal{E}_{h}^{\text{out}}},
\label{eqn:discrete_transport_1} 
\end{multline}
and
\begin{equation}
{\cal{L}}(v) :=
(c_q  (q^{n+1})^+ , v )_{\mathcal{T}_h}
-
\left\langle 
c_{\text{in}} \rho_0  \bU^{n+1} \cdot \bn,  v \right\rangle_{\mathcal{E}_{h}^{\text{in}}}.
\end{equation}
Here the upwind value of concentration is defined by 
$$
{(C^{n+1})}^*_{|e} :=
\left\{
	\begin{array}{ll}
		{(C^{n+1})}^+ & {\mbox{\rm if }}  {\bU^{n+1} \cdot \bn^+ < 0}, \\
		{(C^{n+1})}^- &  {\mbox{\rm if }}  {\bU^{n+1} \cdot \bn^+ \geq 0 },
	\end{array}
\right.
$$
where ${(C^{n+1})}^+$ denotes the value of a neighbor upwind element.
In addition, the source/sink term $q$ for $\tilde{q}=\tilde{c}q$ splits by 
$$
(q^{n+1})^+ = \max(0,q^{n+1}) \quad \text{ and } \quad 
(q^{n+1})^- = \min(0,q^{n+1}),
$$
where $q^{n+1} = (q^{n+1})^+ + (q^{n+1})^-$. 
Recall that $\tilde{c}$ is the injected concentration $c_q$ if $q>0$ and is the resident concentration $c$, if $q<0$.
The $\alpha_c(\tilde{k})$ is {a} penalty parameter for the transport system that can vary on edges
and the choice of $\theta$ leads to different EG algorithms.
For this paper, we also set $\theta = 0$ for the transport problem.  



\subsubsection{Entropy residual stabilization}
\label{sec:levelset_entropy}
The high order transport system ($\tilde{k}\geq 1$) is required to be stabilized to eliminate spurious oscillations 
due to sharp gradients in the exact solution.
In this section, we describe an entropy viscosity  stabilization technique to avoid those oscillations for the EG formulation \eqref{eqn:discrete_transport}.
This method was introduced in \cite{guermond2011entropy} and mathematical stability properties are discussed in 
\cite{MR3167449} for CG and in \cite{Zingan:2013bb} for DG.
First, we start by introducing a numerical dissipation by adding 
\begin{multline}
{\mathcal{E}}(C^{n+1}, v)  
:= 
\left({\mu}^{n+1}_{\text{Stab}}(C,\bU) _{|T} \nabla C^{n+1}, \nabla v \right)_{\mathcal{T}_h}  \\
-
\left\langle \av { {\mu}^{n+1}_{\text{Stab}}(C,\bU) _{|T}  \nabla C^{n+1} }, \jump{v} \right\rangle_{\mathcal{E}_h^{I} }
+
\dfrac{\alpha_s}{h_{e}}\av{ {\mu}^{n+1}_{\text{Stab}}(C,\bU) _{|T}}  \left\langle  \jump{C^{n+1}}, \jump{v} \right\rangle_{\mathcal{E}_h^{I}},
\end{multline}
on the left hand side of the system \eqref{eqn:discrete_transport}.
This results to solve 
\begin{equation}
\mathcal{A}(C^{n+1}, v) + 
\mathcal{E}(C^{n+1}, v)
= {\mathcal{L}}(v), \quad \forall v \in V_{h,\tilde{k}}^{\textsf{EG}}(\mathcal{T}_h).
\label{eqn:discrete_transport_ev} 
\end{equation}
Note that the numerical dissipation term can be added explicitly, but it is known that this would require a time step restriction.
Here ${\mu}^{n+1}_{\text{Stab}}(C,\bU) _{|T}: 
\Omega \times [0,\mathbb{T}] \rightarrow \mathbb{R}$ is
the stabilization coefficient defined on each $T \in \mathcal{T}_h$ by 
\begin{equation}\label{eqn:levelset_stab_v}
{\mu}^{n+1}_{\text{Stab}}(C,\bU) _{|T} := 
\min(
{\mu}^{n+1}_{{\textsf{Lin}}}(C,\bU)_{|T} , 
\mu^{n+1}_{{\textsf{Ent}}}(C,\bU)_{|T} ), 
\end{equation}
which is a piecewise constant over the mesh. 
The main idea here is to split the stabilization: 
when $C(\cdot, t)$ is smooth, the entropy viscosity stabilization $\mu^{n+1}_{{\textsf{Ent}}}(C,\bU)_{|T}$ is activated, 
and when $C(\cdot, t)$ is not smooth because of the complex flux, the linear viscosity ${\mu}^{n+1}_{{\textsf{Lin}}}(C,\bU)_{|T}$ is activated.
The first order linear viscosity is defined by, 
\begin{equation}
{\mu}^{n+1}_{{\textsf{Lin}}}(C,\bU)_{|T}  := \lambda_{{\textsf{Lin}}} \|h \bU^{n+1} \|_{L^{\infty}(T)}, 
\quad \forall T \in \mathcal{T}_h ,
\label{eqn:level:fv}
\end{equation}
where  $\lambda_{{\textsf{Lin}}}$ is a positive constant.

Next, we describe the entropy viscosity stabilization. 
Recall that it is known that the scalar-valued conservation equation
\begin{equation}
{\partial_t} c + \nabla \cdot  {\boldsymbol{f}}(c) = \tilde{q}
\end{equation}
may have one weak solution 
in the distribution sense 
satisfying the additional inequality 
\begin{equation}
{\partial_t} E(c) + \nabla \cdot  {\boldsymbol{F}}(c) - E'(c) \tilde{q} \leq 0, 
\end{equation}
for any convex function $E \in \mathbb{C}^0(\Omega;\mathbb{R})$ which is called entropy and 
${\boldsymbol{F}}'(c):= E'(c) {\boldsymbol{f}}'(c)$, ths associated entropy flux \cite{Kruzkov,Panov1994}. 
The equality holds for smooth solutions. 
We consider ${\boldsymbol{f}}(c) := \bu c$ and obtain 
${\boldsymbol{F}}'(c) = \bu \cdot E'(c)$ and 
$\nabla \cdot {\boldsymbol{F}}(c) = {\boldsymbol{F}}' \cdot \nabla c$.
Note that we can rewrite $\nabla E(C) = E'(C)\nabla C$.
We define the entropy residual which is a reliable indicator of the regularity of $C$ as 
\begin{equation}
R_{\textsf{Ent}}^{n+1}(C,\bU) := 
\bdf{E((C^n)^{\star})}{m} + {\bU^{n+1}} E'((C^n)^{\star})\nabla (C^n)^{\star} - E'((C^n)^{\star}) \tilde{q},
\label{eqn:ent_res}
\end{equation}
which is large when $C$ is not smooth. 
Here $(C^n)^{\star}$ is the extrapolated value which can be utilized as 
$$(C^n)^{\star} :=
\begin{cases}
C^n,     \mbox{ time lagged}, \\
C^n + ( C^n - C^{n-1}), \mbox{ extrapolation}.
\end{cases}
$$
The well known entropy functions include
\begin{equation}
\begin{cases}
E((C^n)^\star) = |(C^n)^\star - r |, \ r \in \mathbb{R}, \text{ Kru\u{z}kov pairs} \\
E((C^n)^\star) = \dfrac{1}{b} |(C^n)^\star|^b, \ $b$  \text{ is a positive even number}  \\
E((C^n)^\star) = - \log( | (C^n)^\star (1-(C^n)^\star)| + \varepsilon), \ \varepsilon \ll 1.
\end{cases}
\label{eqn:entropy_func}
\end{equation}
We chose and test latter two functions in our numerical examples.
Finally, the local entropy viscosity for each step is defined as
\begin{equation}
\mu^{n+1}_{{\textsf{Ent}}}(C, \bU)_{|T} := \lambda_{\textsf{Ent}} h^2 
\dfrac{{ER_{\textsf{Ent}}^{n+1}}_{|T} } 
{\|E((C^n)^\star) - (\bar{E}^n)^\star\|_{L^{\infty}(\Omega)}}, 
\quad  \forall T \in \mathcal{T}_h,
\end{equation}
where 
\begin{equation}
{ER_{\textsf{Ent}}^{n+1}}_{|T} := 
\max(\|R_{\textsf{Ent}}^{n+1}\|_{L^{\infty}(T)}, \|J_{\textsf{Ent}}^{n+1}\|_{L^{\infty}(\partial T \setminus \partial \Omega)}).
\label{eqn:ent_res}
\end{equation}
Here $\lambda_{\textsf{Ent}}$ is a positive constant to be chosen with the average $(\bar{E}^n)^\star := \frac{1}{|\Omega|} \int_{\Omega}E((C^n)^\star) \ d\bx$. 
We define the residual term calculated on the faces by
\begin{equation}
J_{\textsf{Ent}}^{n+1}(C,\bU) := 
h^{-1}_{e}  \av{\bU^{n+1}} \jump{E((C^n)^{\star})}.
\end{equation}
The entropy stability with above residuals for discontinuous case is given 
with  more details
in \cite{Zingan:2013bb}. 
Also, readers refer \cite{guermond2011entropy,JL2016} for tuning the constants 
($\lambda_{\textsf{Ent}}, \lambda_{\textsf{Lin}}$).

\subsection{Adaptive Mesh Refinement}

In this section, we propose a refinement strategy by increasing the mesh resolution in the cells where the entropy residual values \eqref{eqn:ent_res} are higher than others. It is shown in \cite{andrews1998posteriori,puppo2004numerical} that the entropy residual can be used as a posteriori error indicator.
We note that the general residual 
of equation 
\eqref{eqn:discrete_transport_1} 
could also be utilized as an error indicator, but consistency requires that this residual goes to zero as $h \rightarrow 0$.
However, as discussed in \cite{guermond2011entropy}, 
the entropy residual \eqref{eqn:ent_res} converges to a Dirac measure supported in the neighborhood of shocks. 
In this sense, the entropy residual is a robust indicator 
and also efficient since it is been computed for a stabilization.

\begin{figure}[!h]
\centering
\includegraphics[scale=1.3]{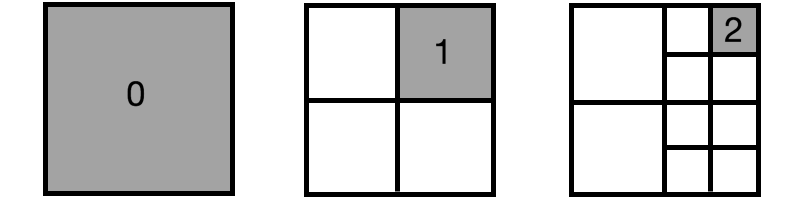}
\caption{Refinement level $\textsf{Ref}_T=$ $0$, $1$, and $2$. 
The number in each cell denotes the each refinement level.}
\label{fig:mesh}
\end{figure}

We denote the number of times a cell($T$) from the initial subdivision has been refined to produce the current cell, 
the refinement level,
by $\textsf{Ref}_T$
(see Figure \ref{fig:mesh}).
Here, a cell $T$ is refined if its corresponding $\textsf{Ref}_T$ is smaller than a given number $R_{\max}$
 and if 
\begin{equation}
|{ER_{\textsf{Ent}}^{n+1}}_{|T}(\bx_T,t)|  \geq C_R,
\end{equation}
where $\bx_T$ is the barycenter of $T$ and $C_R$ is an absolute constant. The purpose of the parameter $R_{\max}$ is to control the total number of cells, which is set to be two more than the initial 
$\textsf{Ref}_T$.
Here $C_R$ is chosen
to mark and refine the cells which represent top 20$\%$ of the values  \eqref{eqn:ent_res} over the domain.
A cell $T$ is coarsen if 
\begin{equation}
|{ER_{\textsf{Ent}}^{n+1}}_{|T}(\bx_T,t)| \leq C_C,
\end{equation}
where
$C_C$ indicates that the cells in the bottom $10\%$ of the values \eqref{eqn:ent_res} over the domain.
However, cell is not coarsen more if the 
$\textsf{Ref}_T$ is smaller than a given number $R_{\min}$.
Here $R_{\min}$ is set to be two less than the initial $\textsf{Ref}_T$.
In addition, a cell is not refine more if the total number of cells are more than $\textsf{Cell}_{\max}$.
{ 
The subdivisions are accomplished with at most one hanging node per face.
During mesh refinement, coefficients are obtained by standard interpolations and in coarsening by restrictions. 
We refer to the documentation of the deal.II library \cite{dealII83} and p4est \cite{p4est} for the computational details. 
}

\subsection{Global Algorithm and Solvers}
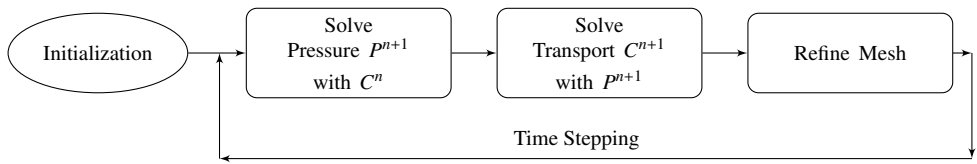
\begin{figure}[H]
\centering
\begin{tikzpicture}[node distance = 3cm, auto]
    \node [cloud] (Init) {\footnotesize Initialization};
    \node [block, right of=Init] (sol_pres) {\footnotesize Solve  \\ Pressure $P^{n+1}$ \\ 
               with $C^{n}$};
    \node [block, right of=sol_pres] (sol_vel) {\footnotesize Solve \\ Transport $C^{n+1}$ \\ 
               with $P^{n+1}$ };
    \node [block, right of=sol_vel] (end) { \footnotesize Refine Mesh};

    \draw[] (6.3,-0.9) node[below]  {\footnotesize  Time Stepping};
    
    \path [line] (Init) -- (sol_pres);
    \path [line] (sol_pres) -- (sol_vel);
    \path [line] (sol_vel) -- (end);

     \draw [line] (1.6, -1.4) -- (1.6, 0);  
     \path [line] (11.25, 0)-- (11.5, 0);    
     \path [line] (11.5,0.) -- (11.5,-1.4); 
     \path [line] (11.5,-1.4) -- (1.6, -1.4);


\end{tikzpicture}
\caption{Global algorithm flowchart including the mesh refinement step.}
\label{fig:arg_flow}
\end{figure}

Here we present our global algorithm in Figure \ref{fig:arg_flow} for modeling the miscible displacement problem.  
The system is solved by an approach formulated in \cite{ewing1983simulation}, where the transport equation is decoupled in time and treated efficiently in a sequential time-stepping scheme.
First, we solve the pressure 
assuming concentration values are obtained by extrapolation of the previous time step values.
Here, we employ the efficient solver which was developed in \cite{LeeLeeWhi15}. Basically, we apply Algebraic Multigrid(AMG) block diagonal preconditioner to the block matrix system with a GMRES solver.
Next, 
we solve the transport equation for concentrations.
The entropy residual which is calculated when solving the transport system is directly employed to refine the mesh. 

\section{Numerical Examples}
\label{sec:num_ex}
This section verifies and evaluates the performance of our proposed algorithm. 
We first demonstrate convergence of the EG transport system with entropy viscosity stabilization in Section 4.1.
The miscible displacement system with dynamic mesh adaptivity 
for solving the coupled flow and the transport in 
two and three dimensional heterogenous porous media is treated in Section 
\ref{subsec:block} - Section \ref{subsec:3d}.
Examples in Section \ref{subsec:fingering} and \ref{sec:ex:radial}  illustrate the effects of viscous fingering instabilities with different viscosity ratios in Hele-Shaw cells with both rectilinear and radial flows.

{ 
The authors developed EG code based on the open-source finite element package deal.II \cite{dealII83} which is coupled with the parallel MPI library \cite{open_mpi} and Trilinos solver \cite{Trilinos-Overview}. 
We employed the dynamic mesh adaptivity feature in deal.II that includes the p4est library \cite{p4est}.
}

\subsection{Example 1. Error Convergence Tests with Single Vortex Problem}
\label{subsec:error}
In this section, we study the convergence of the advection dominated ($d_m =\alpha_l =\alpha_t =0$) EG linear transport system \eqref{eqn:main_transport}-\eqref{eqn:main_bd_def} with the entropy residual stabilization discussed in Section \ref{sec:levelset_entropy}. 
In the unit square $\Omega=(0,1)^2$, the velocity field is given as 
\begin{equation*}
\bu(x,y,t) := \left( 
\begin{array}{r}
-2 \sin(\pi y)\sin(\pi x)^2\cos(\pi y)\cos(\pi t/T_p) \\
2    \sin(\pi x)\sin(\pi y)^2\cos(\pi x)\cos(\pi t/T_p)
\end{array}\right),
\end{equation*}
where the flow is time-periodic which 
$c(\bx,T_p)=c_0(\bx)$, for all $\bx\in \Omega$, 
and we assume $\varphi = \rho_0 =1$.
The initial function $c_0$ is chosen to be the signed distance to the 
circle centered at $(0.5,0.75)$ and of radius $0.15$, i.e
\[
c_0(x,y) := ((x-0.5)^2 + (y-0.75)^2)^\frac12 - 0.15,
\qquad (x,y) \in \Omega.
\]
See Figure~\ref{num:ex1:circel_sv}(a) for the setup. 
The initial $c$ value is transported by the given periodic velocity and reverted to the initial position at the final time $T_p$. This example is referred as a single vortex problem \cite{FLD:FLD4071,lee_thesis}. 
Here, we note that $c$ is a simple tracer value that is transported with a given velocity rather than a physical concentration value.
\begin{figure}[!h]
\centering
\subfloat[n=0]
{\includegraphics[scale=0.18]{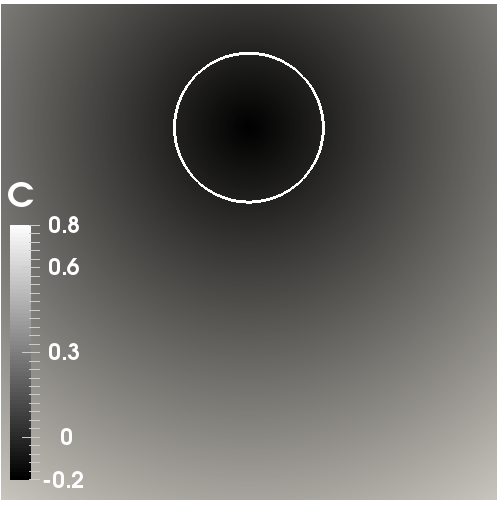}}
\hspace{0.01in}
\subfloat[n=1600]
{\includegraphics[scale=0.18]{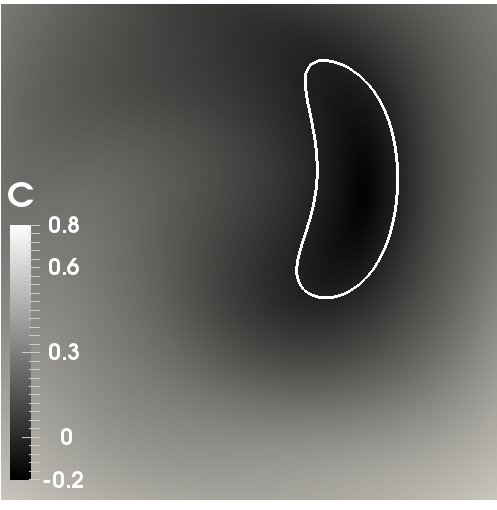}}
\hspace{0.01in}
\subfloat[n=3200]
{\includegraphics[scale=0.18]{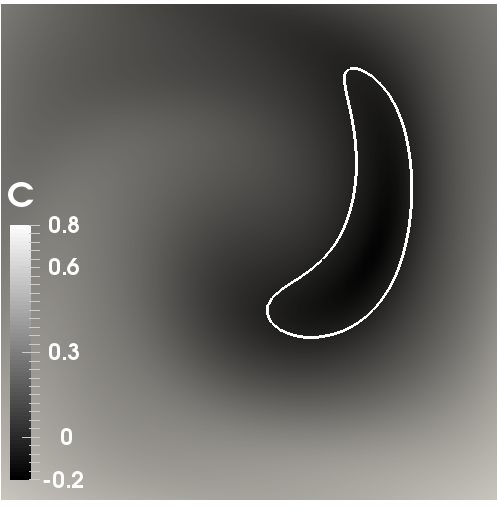}} 
\hspace{0.01in}
\subfloat[n=4800]
{\includegraphics[scale=0.18]{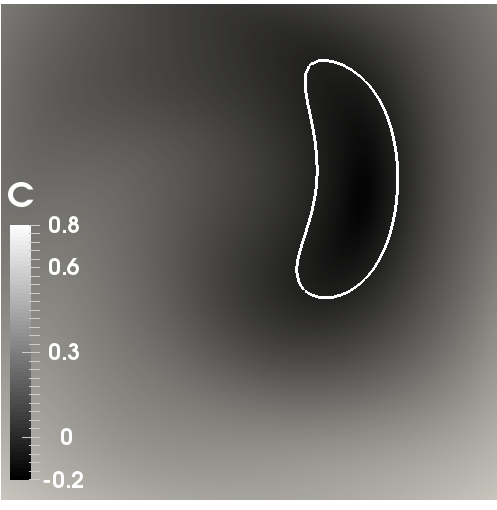}}
\hspace{0.01in}
\subfloat[n=6400]
{\includegraphics[scale=0.18]{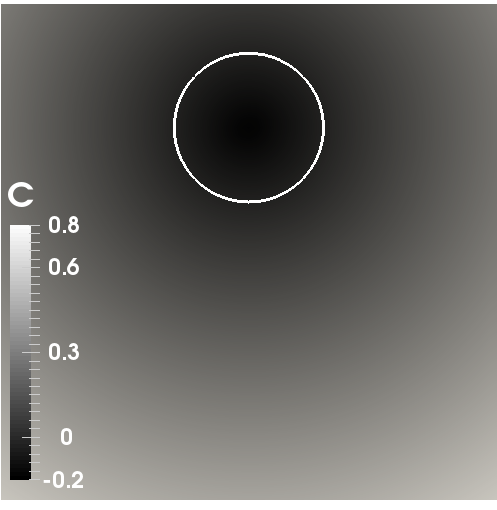}}
\caption{Example 1. Single vortex problem. 
(a)-(e) illustrates the linear transported $c$ values for each time step by a given periodic velocity. The contour line indicates the value $c=0$.}
\label{num:ex1:circel_sv}
\end{figure}
\begin{figure}[!h]
\centering
\subfloat[]
{\includegraphics[width=2.2in]{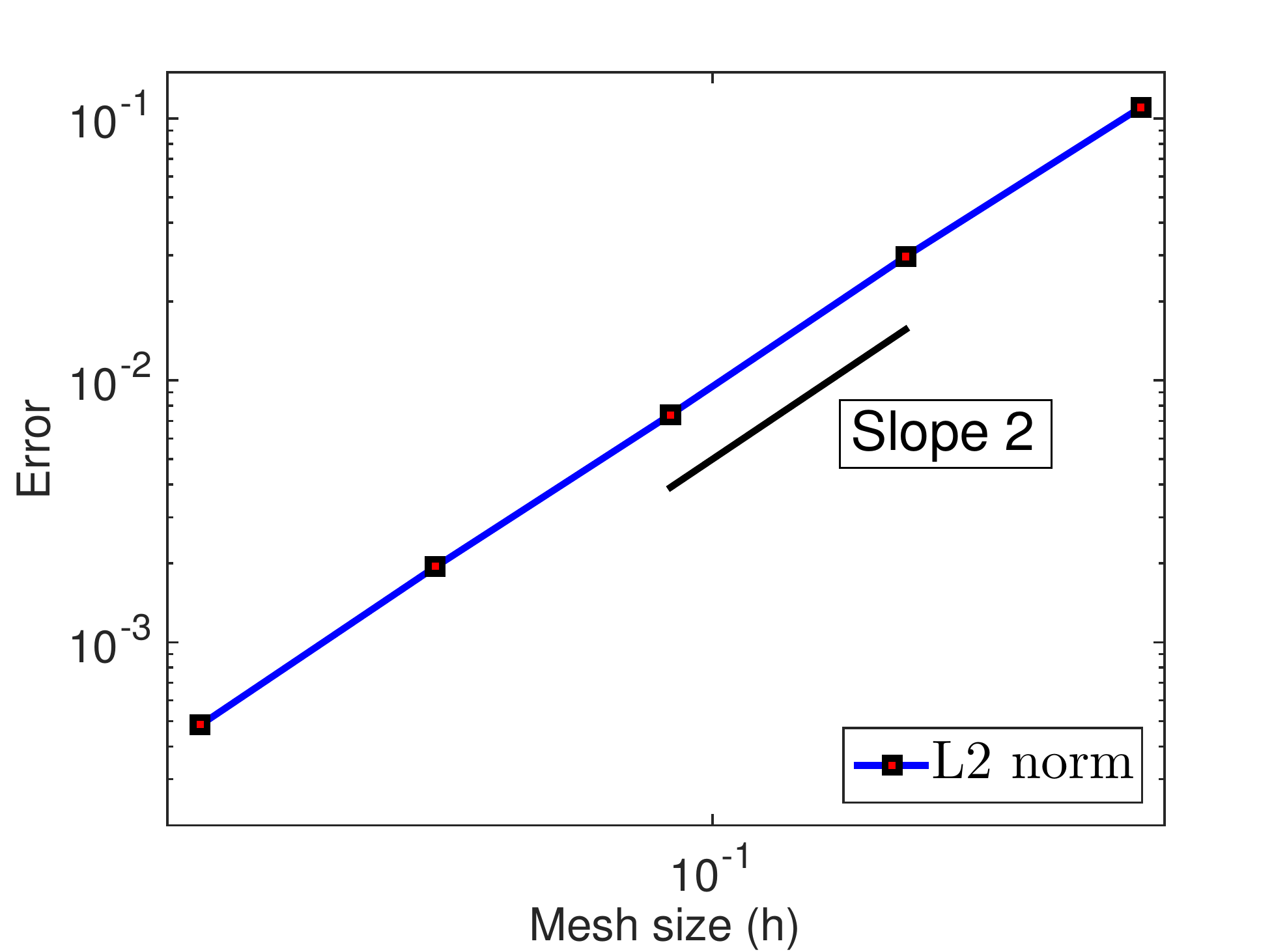}}
\subfloat[]
{\includegraphics[width=2.45in]{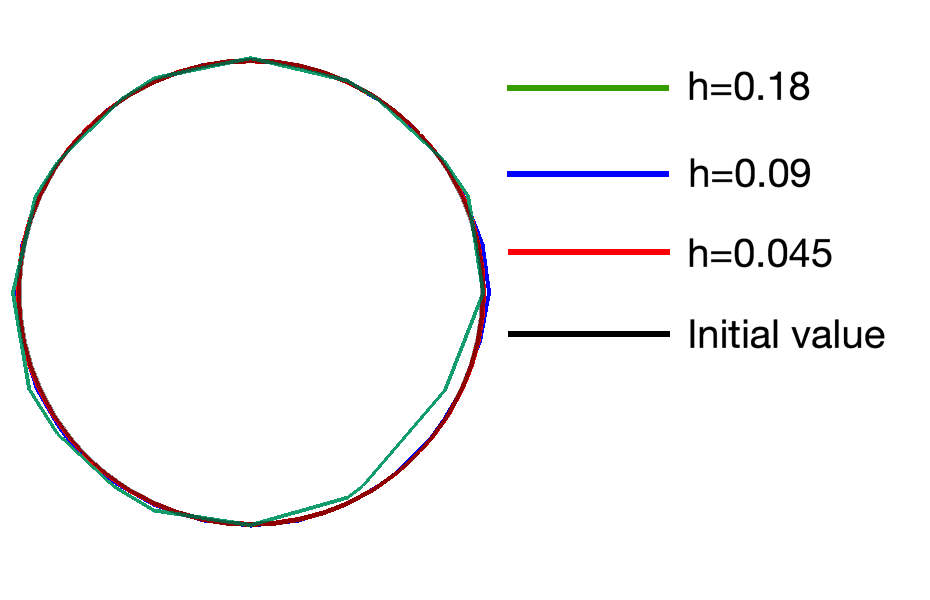}}
\caption{Example 1.Convergence test for the single vortex problem. (a) expected convergence in $\|C^N(\bx, T_p)-C_0\|_{L^2(\Omega)}$ is observed. (b) Comparison of contour values ($c=0$) for each cycle with different $h$ sizes at the final time step; green ($h=0.18$), blue ($h=0.09$), red ($h=0.045$), and black is the initial value. }
\label{fig:ex1:conv}
\end{figure}

We consider five computations on five uniform meshes with constant time steps.  
The mesh-size and the time step are divided by
$2$ each time. 
The meshes are composed of 41, 145, 545, 2113, and 8321
EG-$\polQ_1$ degrees of freedom and 
the time steps are chosen fine
enough not to influence the spacial error.
The errors are evaluated at $t=T_p= 2$.
Here the entropy stabilization coefficients are set to  $\lambda_{\textsf{Lin}} = \lambda_{\textsf{Ent}} = 0.5$ and the entropy function is chosen as 
$E(c) = ({1}/{b}) |c|^b$ with $b=2$.
The expected rates of convergences for the errors 
$\|C^N(\bx, T_p)-C_0\|_{L^2(\Omega)}$
are observed in Figure~\ref{fig:ex1:conv}(a). The comparison of contour values at the final time step for each different $h$ sizes are shown in Figure~\ref{fig:ex1:conv}(b).

\subsection{Example 2. Permeability Block Example using Adaptive Mesh Refinement}
\label{subsec:block}
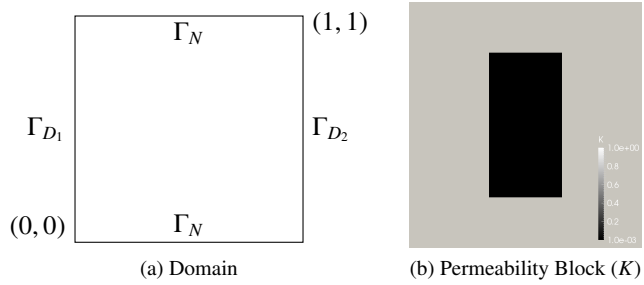
\begin{figure}[!h]
\centering
\subfloat[Domain]{
\begin{tikzpicture}[scale=1.]
\draw[] (0,0) rectangle (3.,3.);
\draw[] (0.,0.2) node[left]     {$(0,0)$};
\draw[] (3.,2.9) node[right]     {$(1,1)$};

\draw[] (1.5,0.5)   node[below]  {$\Gamma_N$};
\draw[] (1.5,2.5)   node[above]   {$\Gamma_N$};

\draw[] (0,  1.5)  node[left]     {$\Gamma_{D_1}$};
\draw[] (3., 1.5)  node[right]  {$\Gamma_{D_2}$};
\end{tikzpicture}
}
\hspace*{0.1in}
\label{fig:domain} 
\subfloat[Permeability Block ($K$)]{\includegraphics[width=1.2in]{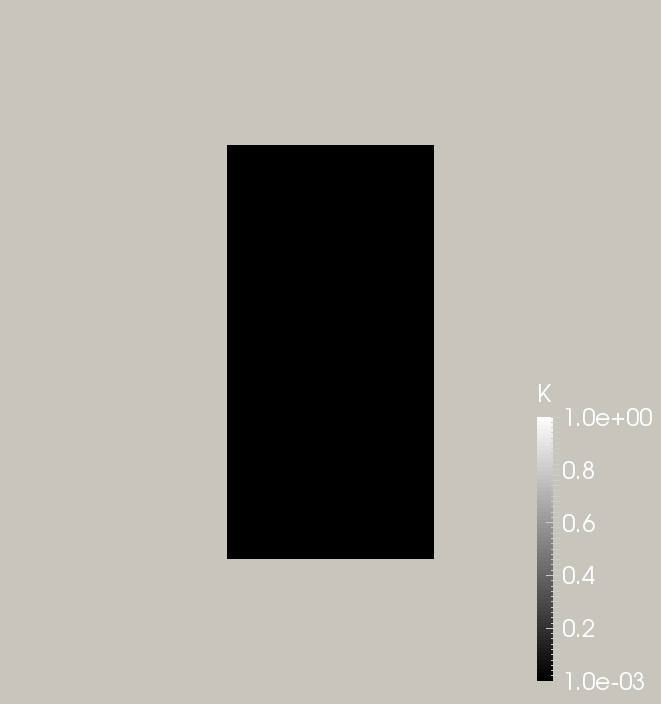}}
\hspace*{0.05in}
\label{fig:k_block}
\caption{Example 2.  (a) computational domain with the boundary conditions 
and (b) permeability block ($K$) in the domain. }
\label{fig:ex1_setup}
\end{figure}
In the computational domain $\Omega = (0,1)^2$, the 
permeability tensor is defined as a diagonal tensor with value $10^{-3}$ in the subdomain 
$\Omega_c = (\frac{3}{8},\frac{5}{8}) \times (\frac{1}{4},\frac{3}{4}) $ and $1$ elsewhere.
Here we set $\varphi = 1$ and 
assume slightly compressible flow by 
$c_F = \num{e-8}$, $\mu = \SI{1}{\pascal \second}$ and $\rho_0 = \SI{1}{\kilogram\per\metre^3}$. See Figure \ref{fig:ex1_setup} for the details with boundary conditions.
We employ EG-$\mathbb{Q}_1$ for the pressure and transport system and we set $d_m =\alpha_l =\alpha_t =0$ for diffusion and dispersion coefficients for this case.

The inflow boundary condition for the transport system \eqref{eqn:main_transport}-\eqref{main_bd_c_out} is given as
$$
c_{\text{in}} = 1 \text{ in \eqref{main_bd_c_in} on }  \Gamma_{D_1} \times (0,\mathbb{T}]  \quad \text{ and } 
\quad 
(-\bD(\bu)\nabla c) \cdot \bn = 0 
\text{ on }\Gamma_N \cup  \Gamma_{D_2}  \times (0,\mathbb{T}], 
$$
and
the initial conditions for the pressure and concentration are set to zero, i.e
$c_0 = 0$ and  $p_0  =0$.
The boundary conditions for the pressure system are chosen as 
$$
p = 1 \ \text{on } \Gamma_{D_1} \times (0,\mathbb{T}], \quad
p = 0 \ \text{on } \Gamma_{D_2} \times (0,\mathbb{T}], \quad 
\text{and } \ 
\bu \cdot \bn = 0 \ \text{on }\Gamma_N \times (0,\mathbb{T}] .
$$
\begin{figure}[!h]
\centering
\subfloat[n=50]
{\includegraphics[width=0.225\textwidth]{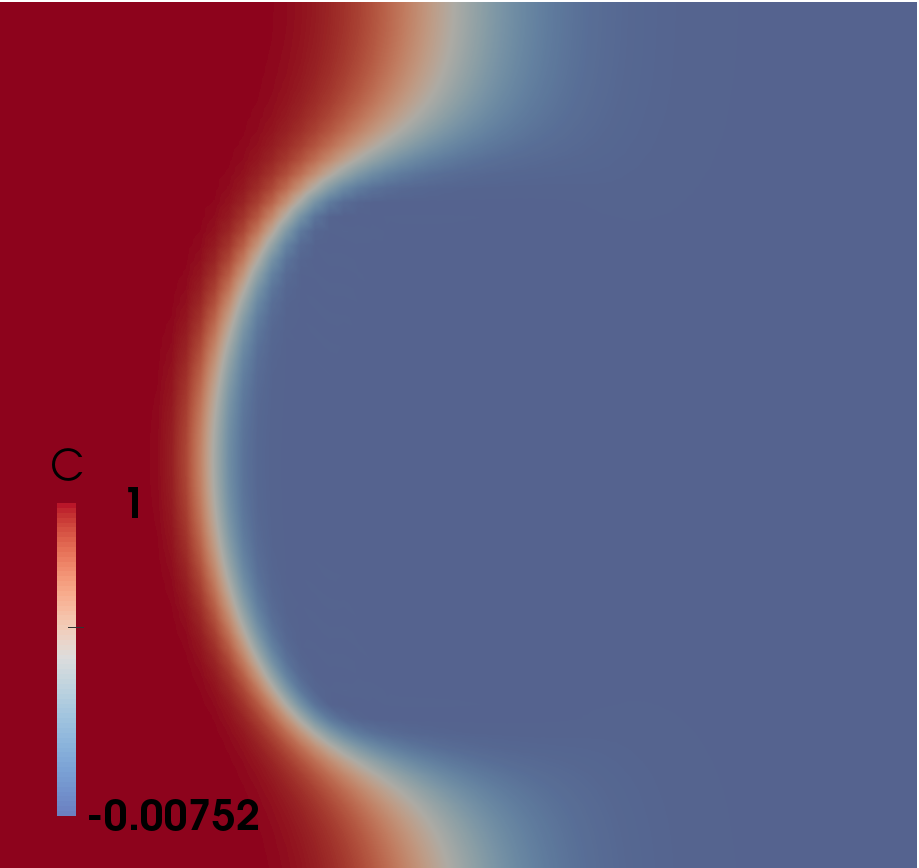}}
\hspace*{0.05in}
\subfloat[n=100]
{\includegraphics[width=0.225\textwidth]{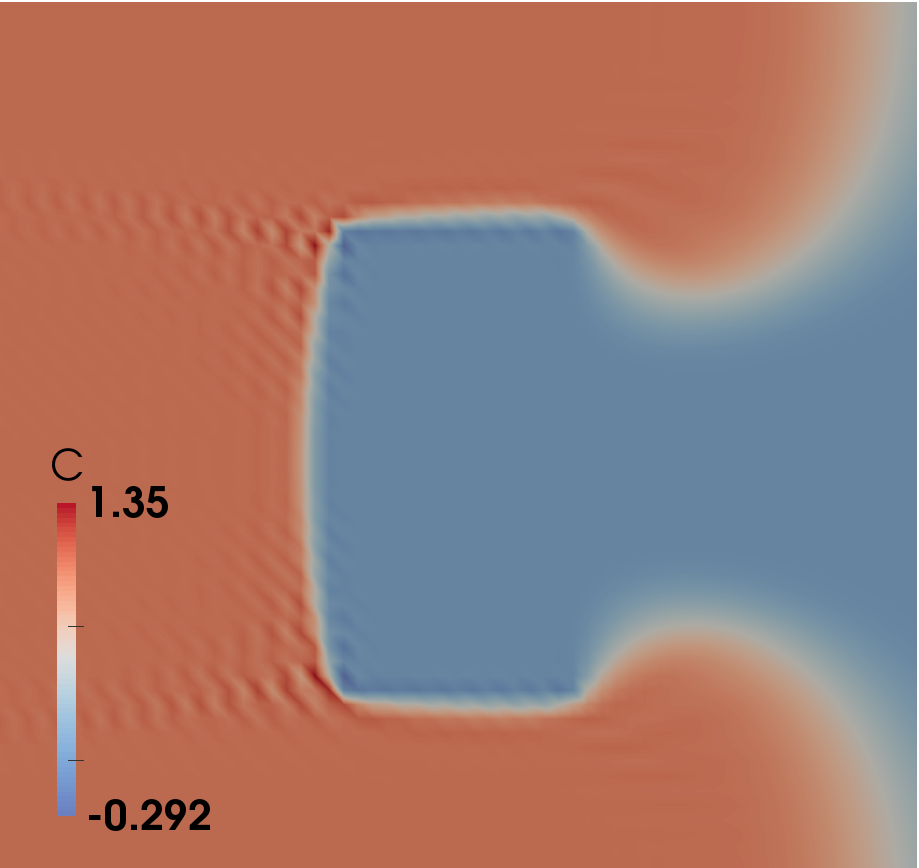}}
\hspace*{0.05in}
\subfloat[n=150]
{\includegraphics[width=0.225\textwidth]{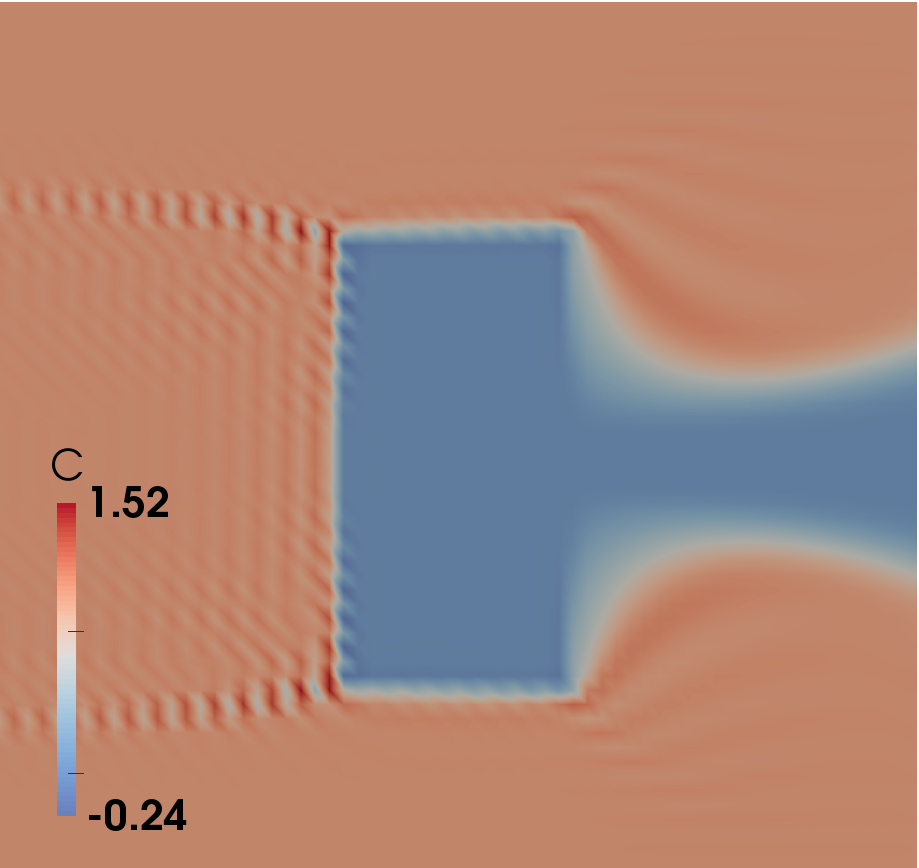}}
\hspace*{0.05in}
\subfloat[n=200]
{\includegraphics[width=0.225\textwidth]{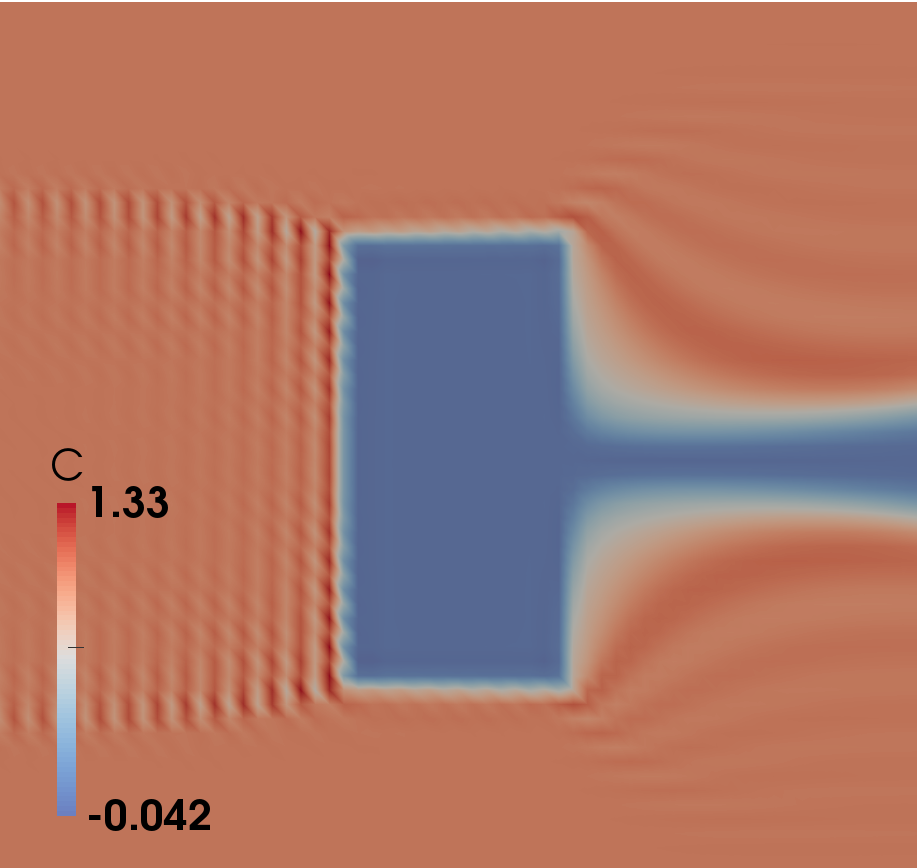}} 
\caption{Example 2. (a)-(d) Concentration values for the each time step $n$ without any entropy residual stabilization, i.e $\lambda_{\textsf{Lin}} = \lambda_{\textsf{Ent}} = 0.$. We observe some spurious oscillations. }
\label{num:ex2:block_no_EV}
\end{figure}
We note that the higher order ($\tilde{k} \geq 1$)  discretization for the transport system requires additional numerical stabilization to avoid any numerical oscillations. 
To demonstrate the impact of the entropy viscosity stabilization,
Figure \ref{num:ex2:block_no_EV} 
illustrates the 
concentration values for each time step without any stabilization.
\begin{figure}[!h]
\centering
\subfloat[n=50]
{\includegraphics[width=0.225\textwidth]{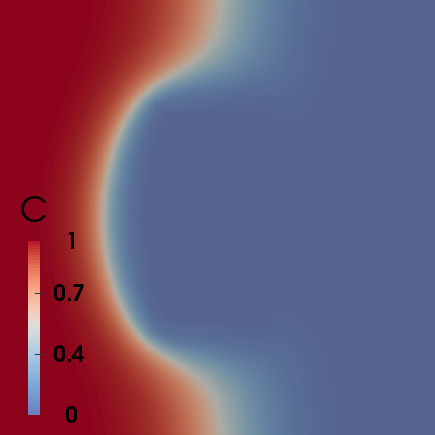}}
\hspace*{0.05in}
\subfloat[n=100]
{\includegraphics[width=0.225\textwidth]{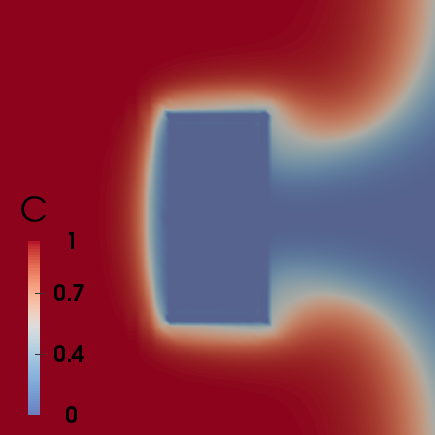}}
\hspace*{0.05in}
\subfloat[n=150]
{\includegraphics[width=0.225\textwidth]{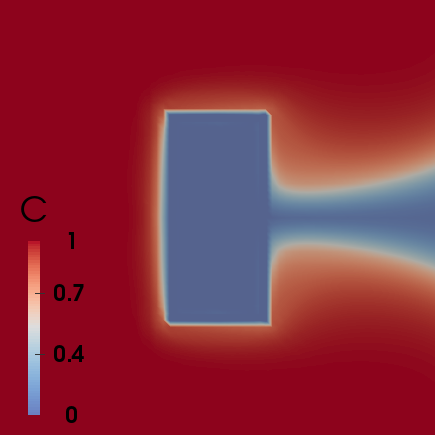}}
\hspace*{0.05in}
\subfloat[n=200]
{\includegraphics[width=0.225\textwidth]{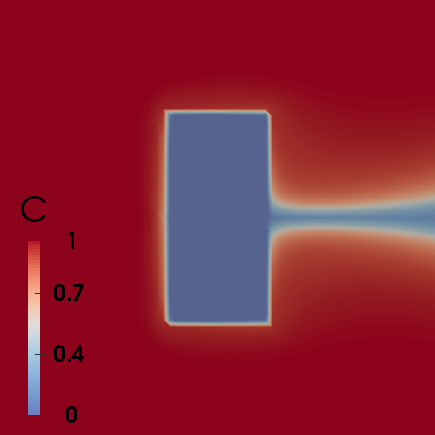}} \\
\subfloat[n=50]
{\includegraphics[width=0.225\textwidth]{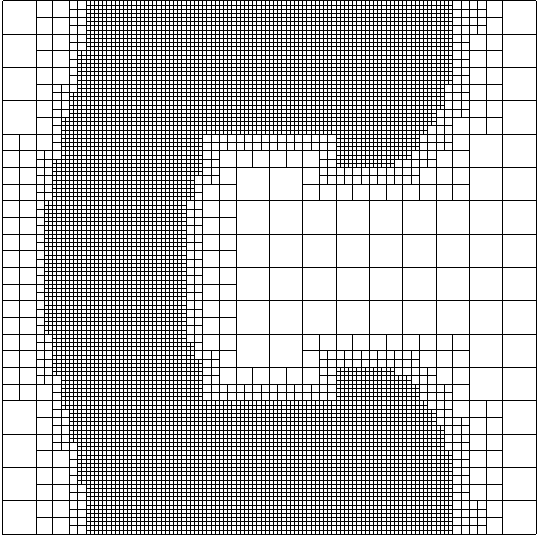}}
\hspace*{0.05in}
\subfloat[n=100]
{\includegraphics[width=0.225\textwidth]{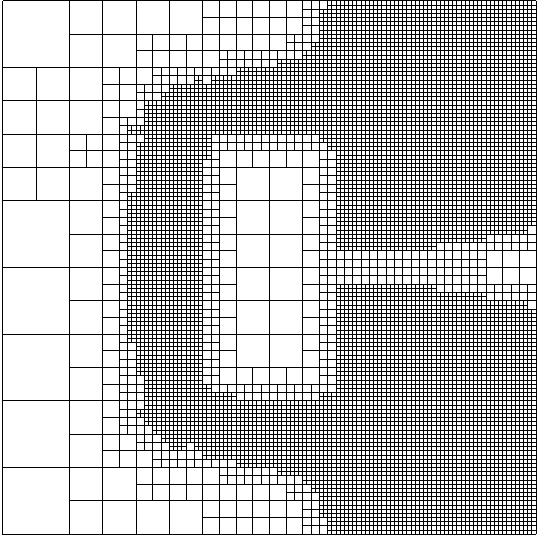}}
\hspace*{0.05in}
\subfloat[n=150]
{\includegraphics[width=0.225\textwidth]{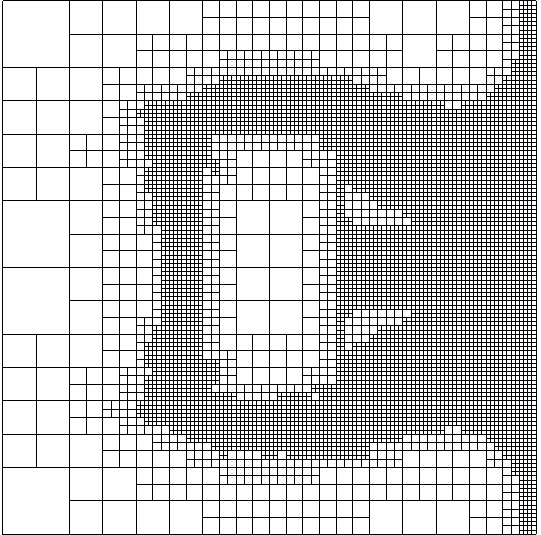}}
\hspace*{0.05in}
\subfloat[n=200]
{\includegraphics[width=0.225\textwidth]{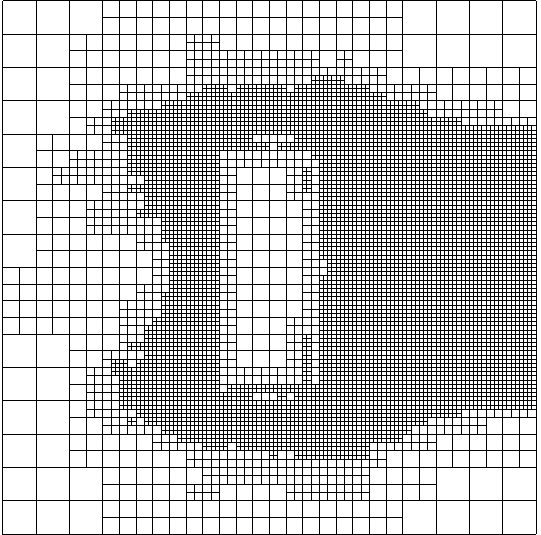}}
\caption{Example 2. 
(a)-(d) concentration values for the each time step $n$. Oscillations are avoided by applying the entropy residual stabilization, $\lambda_{\textsf{Lin}} = \lambda_{\textsf{Ent}} = 0.5$. 
(e)-(h) corresponding adaptively refined meshes for each time step are also illustrated at the bottom row.}
\label{num:ex2:block_EV}
\end{figure}

\begin{figure}[!h]
\centering
\subfloat[]
{\includegraphics[width=0.225\textwidth]{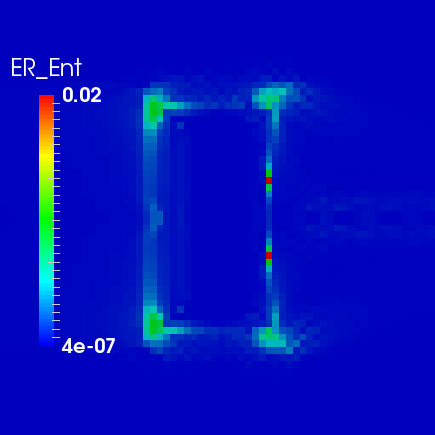}}
\hspace*{0.05in}
\subfloat[]
{\includegraphics[width=0.225\textwidth]{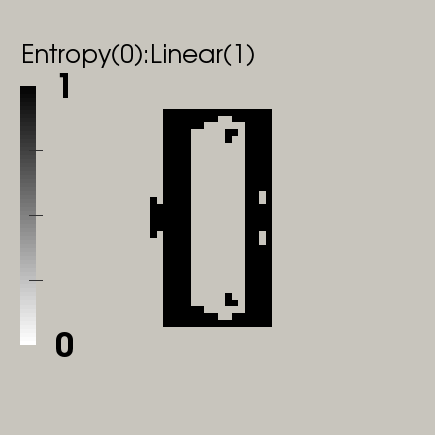}}
\hspace*{0.05in}
\caption{Example 2. (a) entropy residual values $ER_{Ent}$ per cell and (b) selection of the stabilization for \eqref{eqn:levelset_stab_v} at $n=200$. We observe the expected choices for the stabilization. }
\label{num:ex2:block_EV_Entropy}
\end{figure}

Next, we employ the entropy residual stabilization discussed in previous section 
with an 
entropy function chosen as $E(c) = - \log( | c (1-c)| + \varepsilon)$, $\varepsilon= \num{e-4}$. Here the stabilization coefficients are set to 
$\lambda_{\textsf{Lin}} = \lambda_{\textsf{Ent}} = 0.5$. 
The results are shown at Figure \ref{num:ex2:block_EV} without any oscillations.
The numerical discretizations are given as $h_{\min}=0.02$ and $\delta t = 0.01$.
Here $h_{\min}$ denotes the minimum mesh size over the domain with the adaptive mesh refinement.
 Due to the dynamic mesh refinement, the EG degrees of freedom for transport is around $15,000$ with $\textsf{Cell}_{\max} = 7500$,
$\textsf{R}_{\max} = 7$, 
and $\textsf{R}_{\min} = 3$.
In addition, the Figure \ref{num:ex2:block_EV_Entropy} (a) illustrates the values of the entropy residual defined by \eqref{eqn:ent_res} for each cell. As expected, the values are higher near the jumps (or shocks). 
The choice of the stabilization described in  \eqref{eqn:levelset_stab_v} is shown 
at Figure \ref{num:ex2:block_EV_Entropy} (b).
Here the region with zero (white) values 
and one (block) values indicate, 
where the entropy residual stabilization is activated and 
the linear stabilization is activated, respectively. 
We note that the linear viscosity is chosen at the interfaces where we have large entropy values.

\subsection{Example 3. Random Permeability Tensor (2D) }
\label{subsec:2d}

In this example, we consider miscible displacement flow problem in a two dimensional  heterogenous porous media. The computational domain is $\Omega = (\SI{0}{\metre},\SI{1}{\metre})^2$ and the random permeability tensor is given by
\begin{equation}
K(\bx) = \min \left( \max \left( \sum_{i=1}^M \sigma_i(\bx)  ,0.01\right)   , 4 \right),  \quad 
\sigma_i(\bx) = \exp \left(- \left( \dfrac{|\bx - \bx_i|}{0.05} \right)^2 \right),
\label{eqn:random_perm}
\end{equation}
where the centers $\bx_i$ are $M$ randomly chosen locations inside the domain (\cite{deal}, example step-21). 
{ The latter is assigned on adaptive grids directly on the multiple levels.}
In addition, we set the diffusion and dispersion tensor with the physical coefficient chosen as 
\begin{equation}
d_m = 1.8e^{-7} m^2/s,
\alpha_l = 1.8e^{-5} m^2/s, \mbox{ and }
\alpha_t = 1.8e^{-6} m^2/s.
\end{equation}
All the other physical and numerical parameters and boundary conditions are the same as in the previous example. 
Figure \ref{num:ex3:random_2d} illustrates the EG-$\polQ_1$ solution of concentration values for each time step with corresponding mesh refinements.  
Here the stabilization coefficients are set to 
$\lambda_{\textsf{Lin}} = \lambda_{\textsf{Ent}} = 0.5$. 
We observe that the mesh is refined near the sharp interfaces.
The numerical parameters chosen are  
$h_{\min}=0.02$ and $\delta t = 0.01$ 
with 
$\textsf{R}_{\max} = 7$ and 
$\textsf{R}_{\min} = 3$.

\begin{figure}[!h]
\centering
\subfloat[n=100]
{\includegraphics[width=0.21\textwidth]{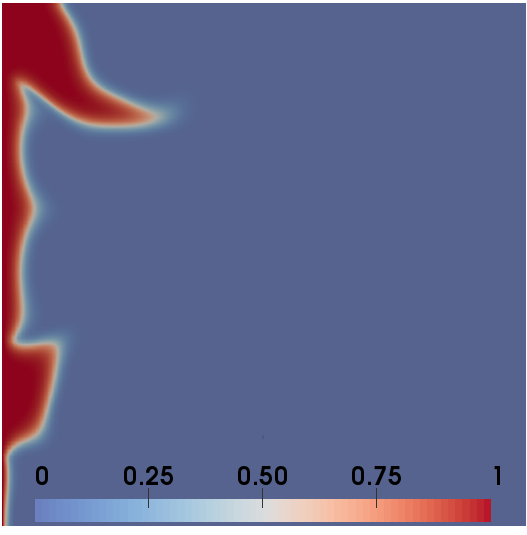}}
\hspace*{0.01in}
\subfloat[n=200]
{\includegraphics[width=0.21\textwidth]{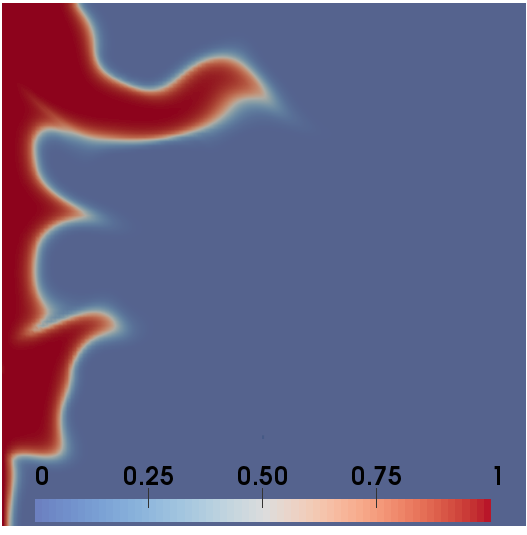}}
\hspace*{0.01in}
\subfloat[n=300]
{\includegraphics[width=0.21\textwidth]{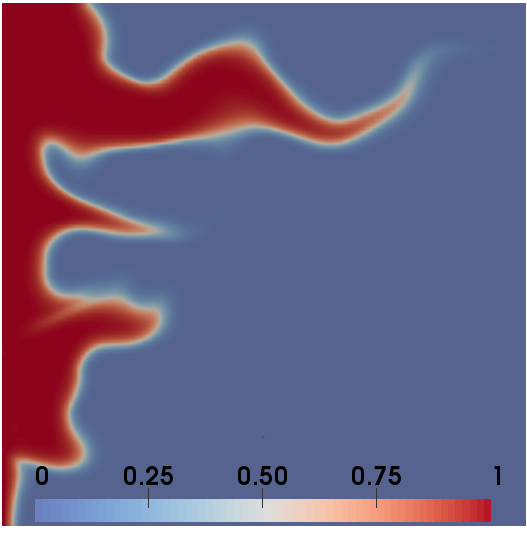}} 
\hspace*{0.01in}
\subfloat[n=400]
{\includegraphics[width=0.21\textwidth]{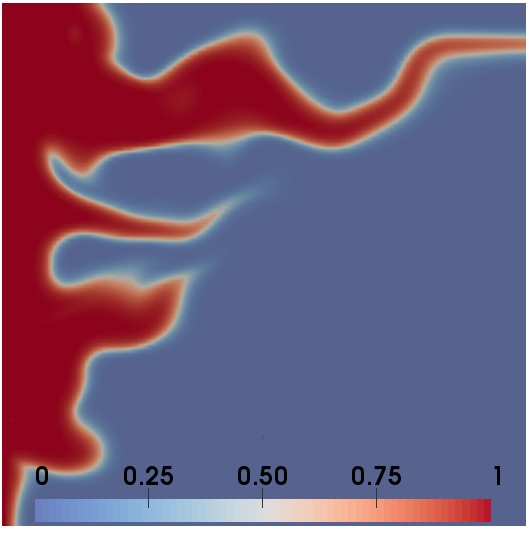}}\\
\subfloat[n=100]
{\includegraphics[width=0.21\textwidth]{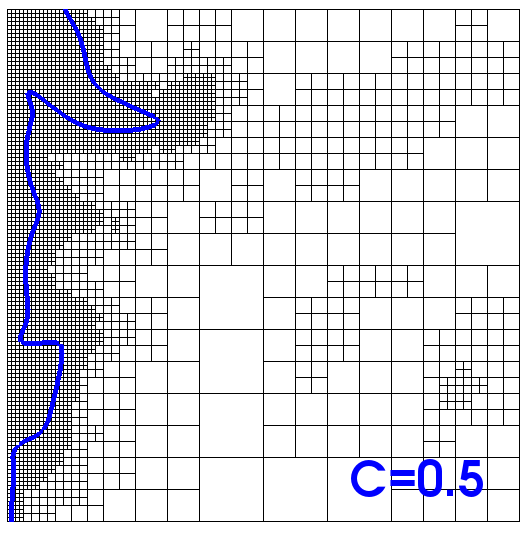}}
\hspace*{0.01in}
\subfloat[n=200]
{\includegraphics[width=0.21\textwidth]{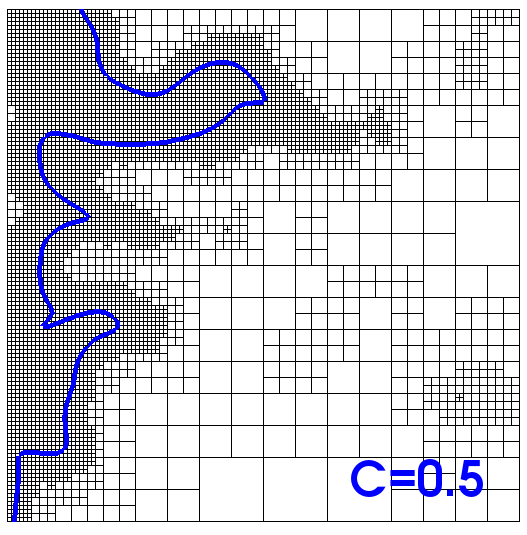}}
\hspace*{0.01in}
\subfloat[n=300]
{\includegraphics[width=0.21\textwidth]{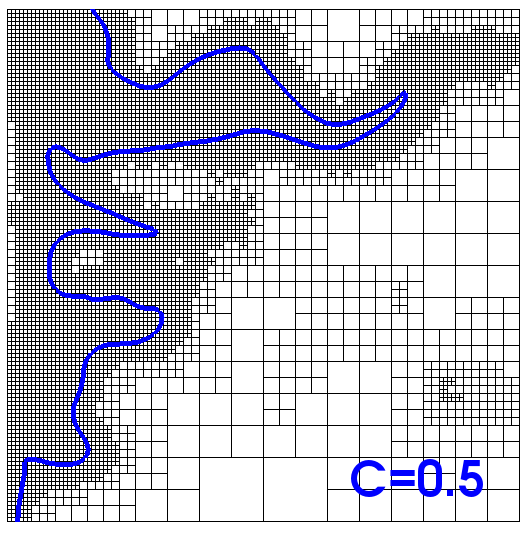}} 
\hspace*{0.01in}
\subfloat[n=400]
{\includegraphics[width=0.21\textwidth]{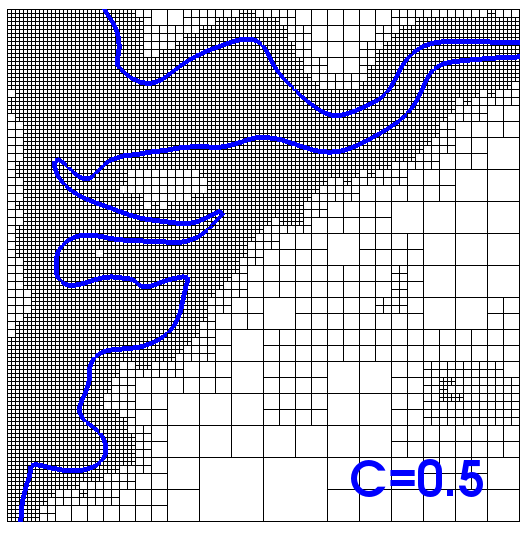}}
\caption{Example 3. (a)-(d) concentration values for each time by a given random tensor permeability. (e)-(h) each corresponding mesh refinements. The blue contour line indicates the value $C=0.5$.}
\label{num:ex3:random_2d}
\end{figure}

\begin{figure}[!h]
\centering
\subfloat[Contour values]
{\includegraphics[width=0.21\textwidth]{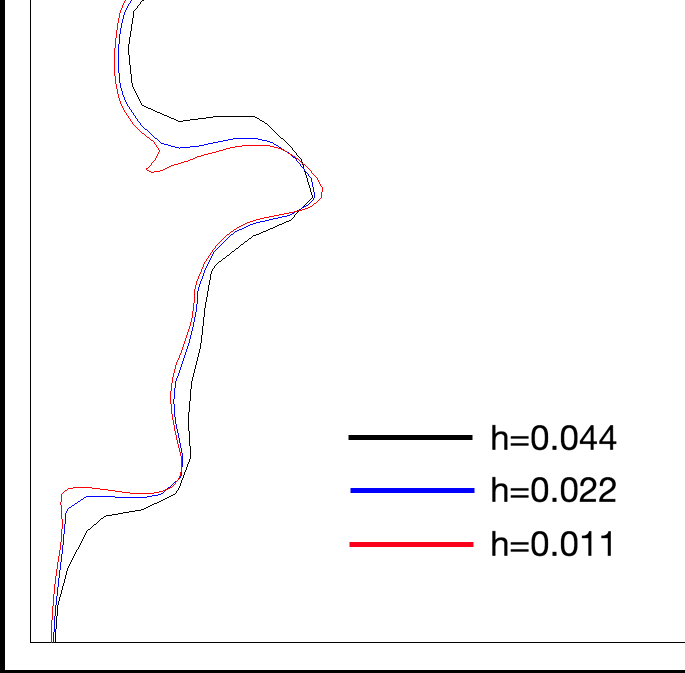}}
\hspace*{0.1in}
\subfloat[Mesh refinement]
{\includegraphics[width=0.21\textwidth]{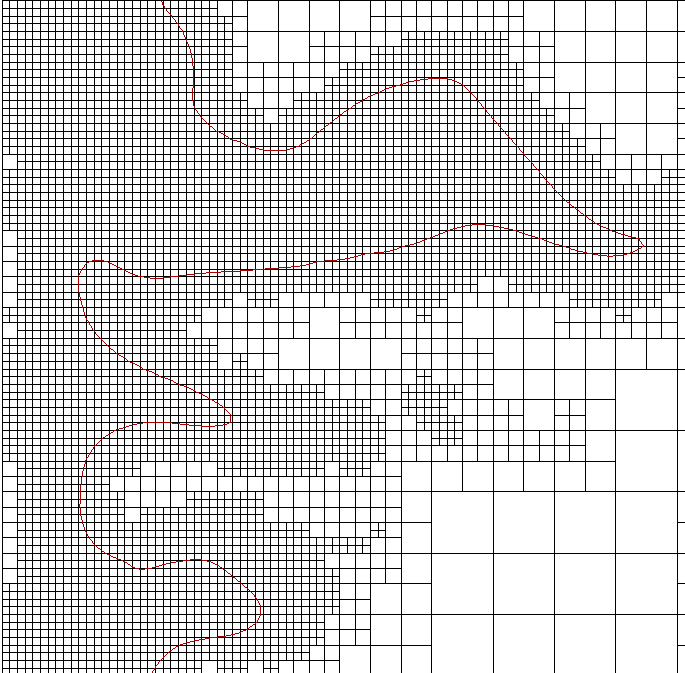}}
\caption{Example 3. 
a) illustrates the contour values ($C=0.5$) for each different mesh sizes ($h=0.044$, $h=0.022$, and $h=0.011$). b) magnified the adaptive mesh refinement for the top left corner at $n=250$. We observe refined and coarsened cells.}
\label{num:ex3:random2d_add}
\end{figure}

In Figure \ref{num:ex3:random2d_add}, we capture the contour of concentration value ($C=0.5$)  at the bottom left part of the domain with different mesh sizes. In addition, we magnify the top left part of the domain to see the mesh refinement at time step $n=250$.

%

\subsection{Example 4. Random Permeability Tensor (3D)}
\label{subsec:3d}
In this section, we consider three dimensional computational domain $\Omega = (\SI{0},\SI{1}{\metre})^3$ and 
the previously defined random permeability tensor \eqref {eqn:random_perm} is multiplied by $2.5$. 
Here the numerical parameters are chosen as $h_{\min}=0.1$ and $\delta t = 0.01$, but 
all the other physical and numerical parameters and boundary conditions are the same as in the previous example. 
Due to the dynamic mesh refinement ($\textsf{R}_{\max} = 6$ and $\textsf{R}_{\min} = 2$), the number of degrees of freedom for EG transport are approximately $240,000$ with the maximum number of cells equal to $112,225$ over the entire time period. 
In particular, three dimensional examples are computed by employing multiple parallel processors (MPI). 
Figure \ref{num:ex3:random_3d} illustrates the EG-$\polQ_1$ solution of concentration values for each time step with corresponding mesh refinements. The adaptive mesh refinement strategy becomes very efficient for large-scale three dimensional problems. 

\begin{figure}[!h]
\centering
\subfloat[$n=4$]
{\includegraphics[width=0.25\textwidth]{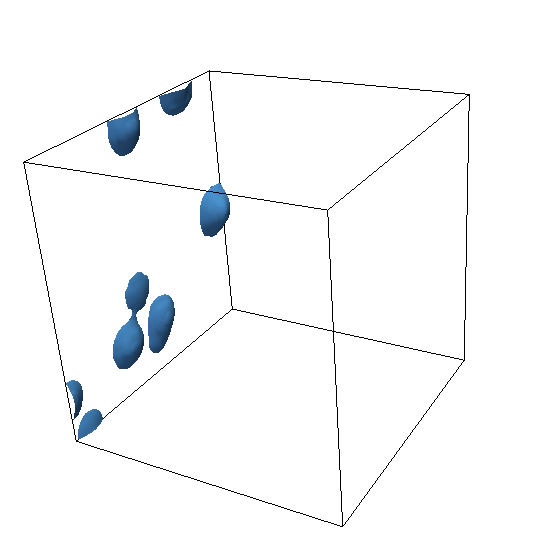}}
\subfloat[$n=100$]
{\includegraphics[width=0.25\textwidth]{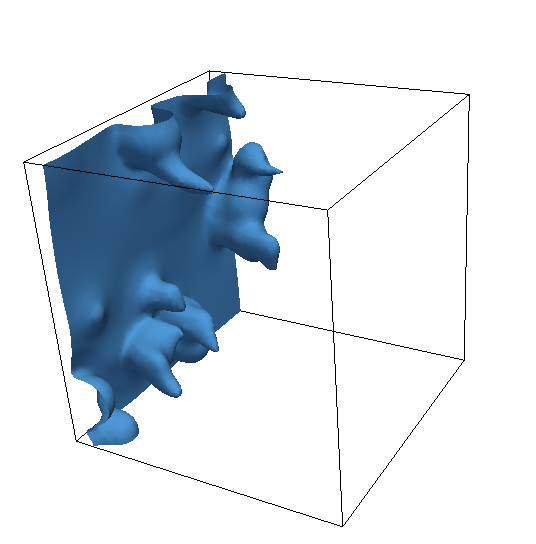}}
\subfloat[$n=160$]
{\includegraphics[width=0.25\textwidth]{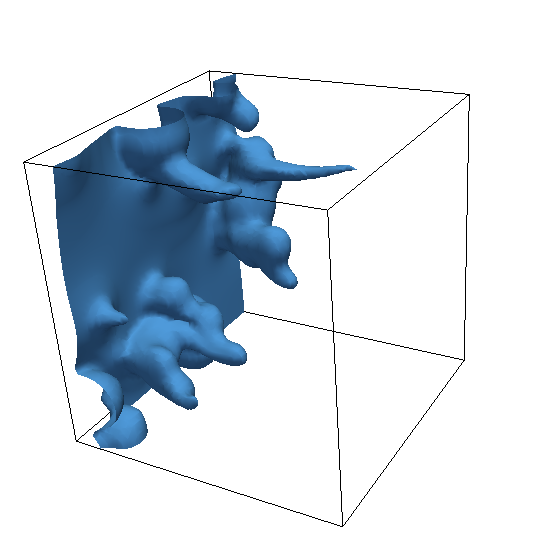}} \\
\subfloat[$n=200$]
{\includegraphics[width=0.25\textwidth]{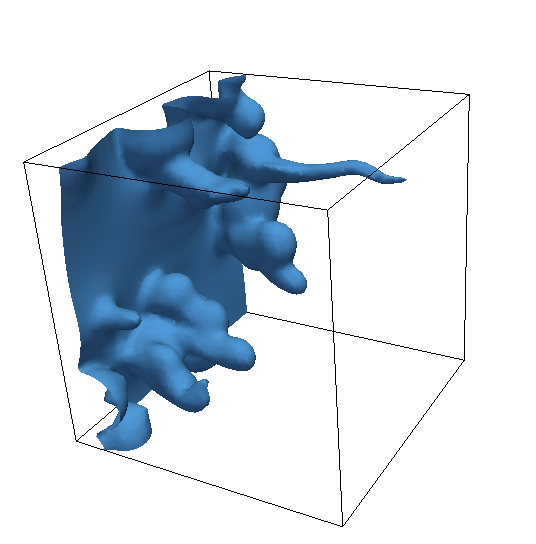}}
\subfloat[$n=285$]
{\includegraphics[width=0.25\textwidth]{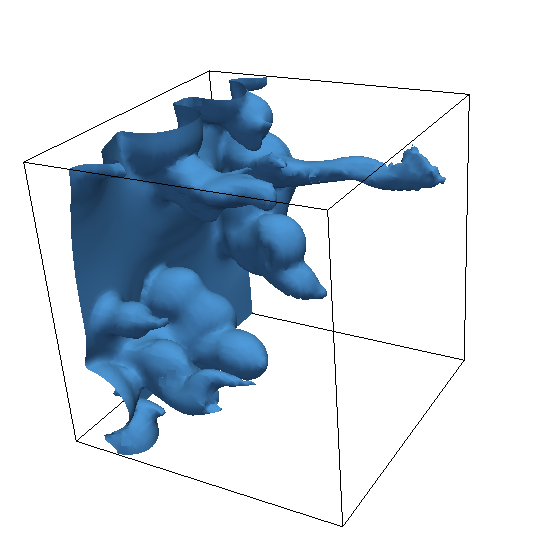}}
\subfloat[$n=285$]
{\includegraphics[width=0.25\textwidth]{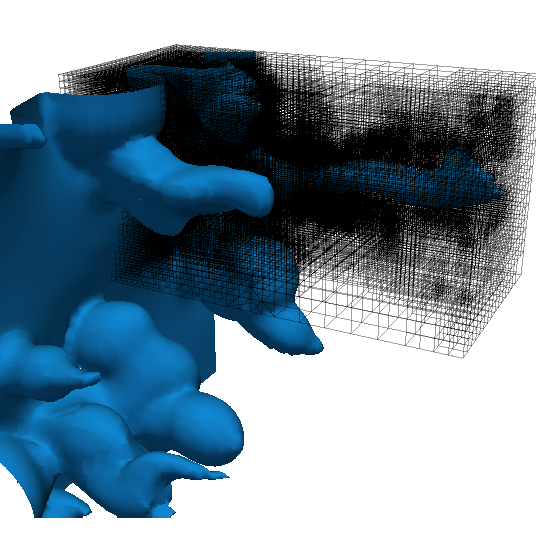}}
\caption{Example 4. The iso-surface indicates the concentration value ($C=0.5$) in  the three dimensional heterogeneous media with a given random tensor permeability. The last figure (f) emphasizes the mesh refinement along the finger.}
\label{num:ex3:random_3d}
\end{figure}


\subsection{Example 5. Hele-Shaw cell: viscous fingering in a homogeneous channel}
\label{subsec:fingering}
\begin{figure}[!h]
\centering
{\includegraphics[scale=0.38]{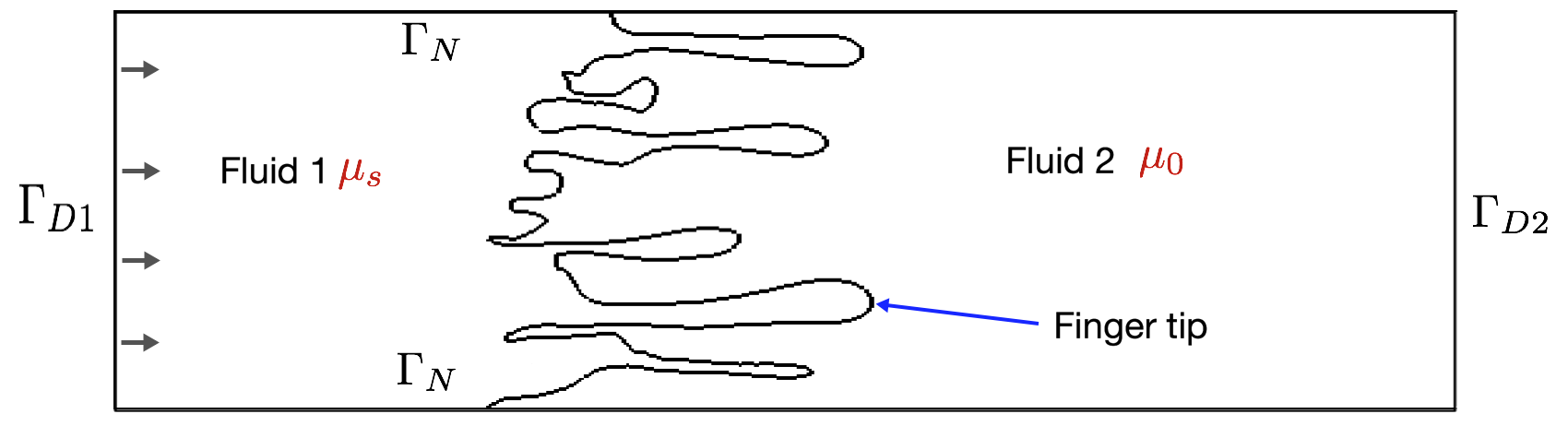}}
\caption{Example 5. Computational domain with  boundary conditions. Mixing zone length and fastest finger tip is defined.}
\label{fig:finger_setup}
\end{figure}
\begin{table}[!h]
\centering
  \begin{tabular}{ c | c | c | c | c | c | c | c | c | c | c }
    \hline
    Case & 1 & 2 & 3 & 4 & 5 & 6 & 7 & 8 & 9 & 10   \\ \hline\hline
  $M$  & 1 & 25 & 50 & 100 & 125 & 150 & 250 & 300 & 500 & 750 \\ \hline
  $\mu_0 \ (\text{Pa s})$  & 0.001 & 0.025 & 0.05  & 0.1 & 0.125 & 0.15 & 0.250 & 0.3 & 0.5 & 0.75 \\ \hline
       \end{tabular}
\caption{Example 5. Different test cases with various viscosity ratios $M$.}
\label{table:vf}
\end{table}

In this example, we present rectilinear flow displacement problems in a Hele-Shaw cell with a different viscosity ratios by injecting water into viscous fluids.
The computational domain $\Omega = (0,0) \times (L,H)$ is defined with 
$(L,H):= (\SI{1}{\metre},\SI{0.25}{\metre})$, 
and
the initial conditions for the pressure and concentration are set to zero, i.e
\begin{equation}
c_0 = 0  \  \text{ and  }  \  p_0  =0.
\end{equation}
The boundary conditions for the pressure and transport system are given as 
\begin{equation}
p = p_{\text{in}} \text{ on } \Gamma_{D_1} \times (0,\mathbb{T}], \quad
p = p_{\text{out}}  = 0 \text{ on } \Gamma_{D_2} \times (0,\mathbb{T}], \quad 
\bu \cdot \bn = 0 
\text{ on }\Gamma_N \times (0,\mathbb{T}], 
\end{equation}
and 
\begin{equation}
c_\text{in} = 1 
\text{ for (9) on } \Gamma_{D_1} \times (0,\mathbb{T}], \ \ 
(-\bD(\bu)\nabla c) \cdot \bn = 0 
\text{ on }\Gamma_N \cup  \Gamma_{D_2}  \times (0,\mathbb{T}], 
\end{equation}
respectively. 
See Figure \ref{fig:finger_setup} for more details. 
In addition, the physical coefficients for the diffusion and the dispersion tensor are chosen as 
\begin{equation}
d_m = 1.8e^{-8} m^2/s,
\alpha_l = 1.8e^{-8} m^2/s, \mbox{ and }
\alpha_t = 1.8e^{-9} m^2/s,
\label{eqn:diff_disp_finger}
\end{equation}
and the fluid is assumed to be incompressible with $c_F = 0$ and $\rho_0 = \SI{1000}{\kilogram\per\metre^3}$.
{ The numerical stabilization coefficients are set to 
$\lambda_{\textsf{Lin}} = 0.8$ and  $\lambda_{\textsf{Ent}} = 0.9$. }
Here $K = I$ and we apply the following 
quarter-power mixing rule \cite{koval1963method} 
for the viscosity 
\begin{equation}
\mu(c) := (c \mu_s^{-0.25} + (1-c)\mu_0^{-0.25})^{-4},
\end{equation}
where $\mu_s$ is the solvent and $\mu_0$ is residing viscosity. Here we denote the viscosity ratio as 
$$
M := \dfrac{\mu_0}{\mu_s}.
$$
Ten different cases of 
viscosity ratios  are
given in Table \ref{table:vf}.
For the inflow boundary condition, we note that $p_{\text{in}}$ is set by 
$$
\dfrac{K}{\mu_0}  \left(\dfrac{p_{\text{in}}-p_{\text{out}}}{L}  \right) = 0.05,
$$
for each case 
in order to obtain a constant initial velocity (0.05) of the residing fluid.
Viscous fingering occurs due to the heterogeneous permeability as shown in the previous examples 
but also can occur for the viscosity ratios $M>1$, which is referred as a Hele-Shaw problem. 
In this latter case the instabilities highly dependent  on $M$ and the P\'{e}clet number which is defined as 
$$
Pe  = \dfrac{L \bU_L}{d_m},
$$
where $L$ is the characteristic length and $\bU_L$ the local flow velocity. 
The numerical parameters for these examples are set to $\delta t= 0.01$ and $h_{\min} = 0.014$ with stability coefficients 
$\lambda_{\textsf{Lin}} = \lambda_{\textsf{Ent}} = 1$. 
\begin{figure}[!h]
\centering
\subfloat[t=5]
{\includegraphics[width=0.45\textwidth]{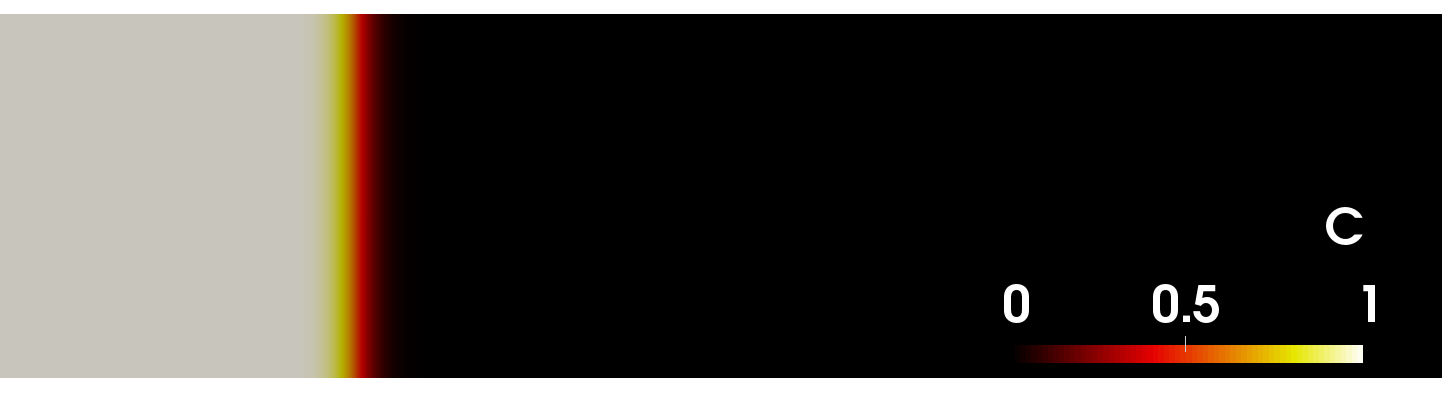}} 
\hspace*{0.01in}
\subfloat[t=15]
{\includegraphics[width=0.45\textwidth]{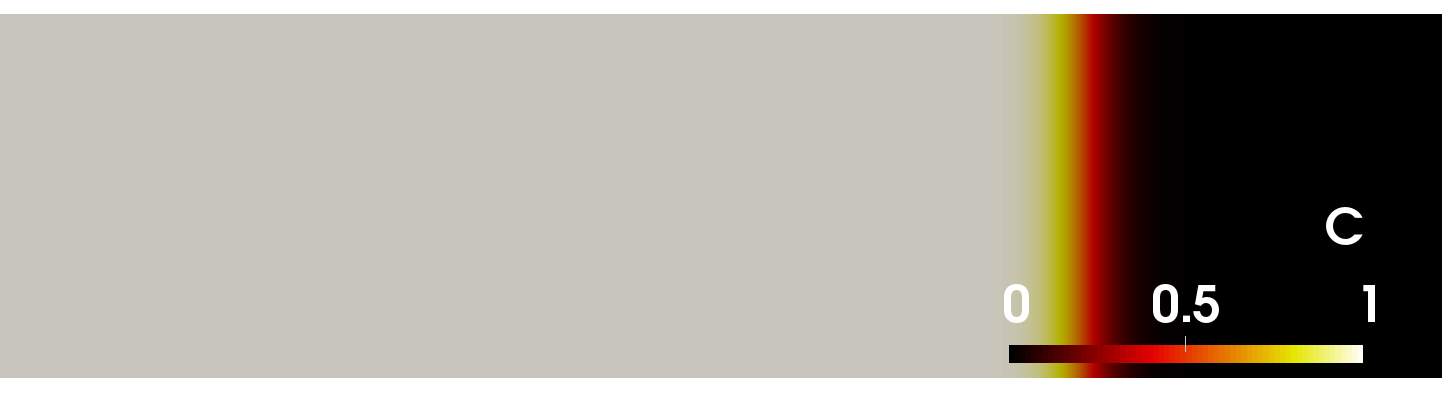}}
\caption{Example 5. Case 1 with $M=1$. Concentration values are shown for each time. The existing fluid is replaced by the same fluid smoothly without any instabilities. }
\label{fig:ratio_1}
\end{figure}
\begin{figure}[!h]
\centering
\subfloat[t=1]
{\includegraphics[width=0.45\textwidth]{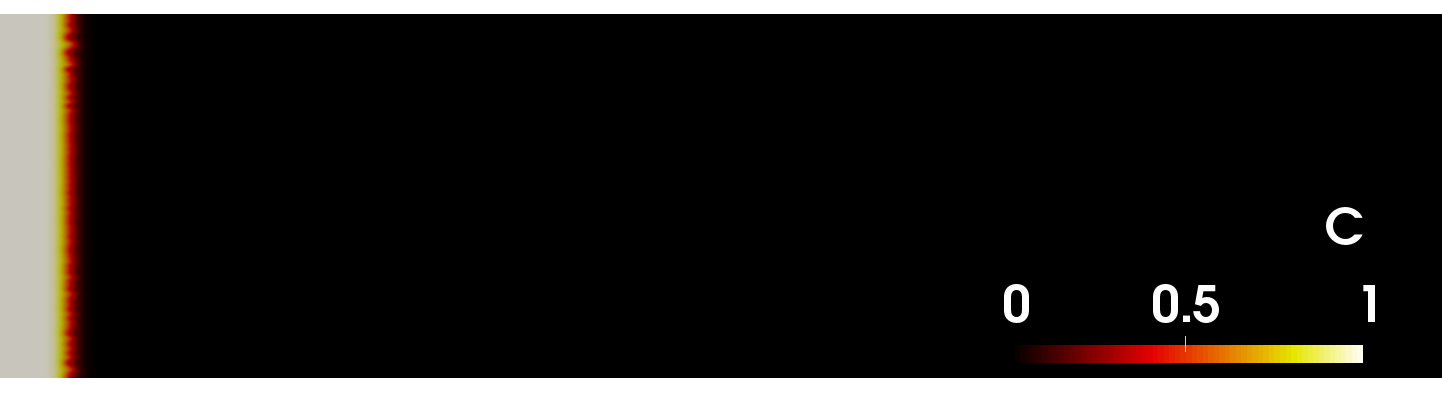}} 
\hspace*{0.01in}
\subfloat[t=2]
{\includegraphics[width=0.45\textwidth]{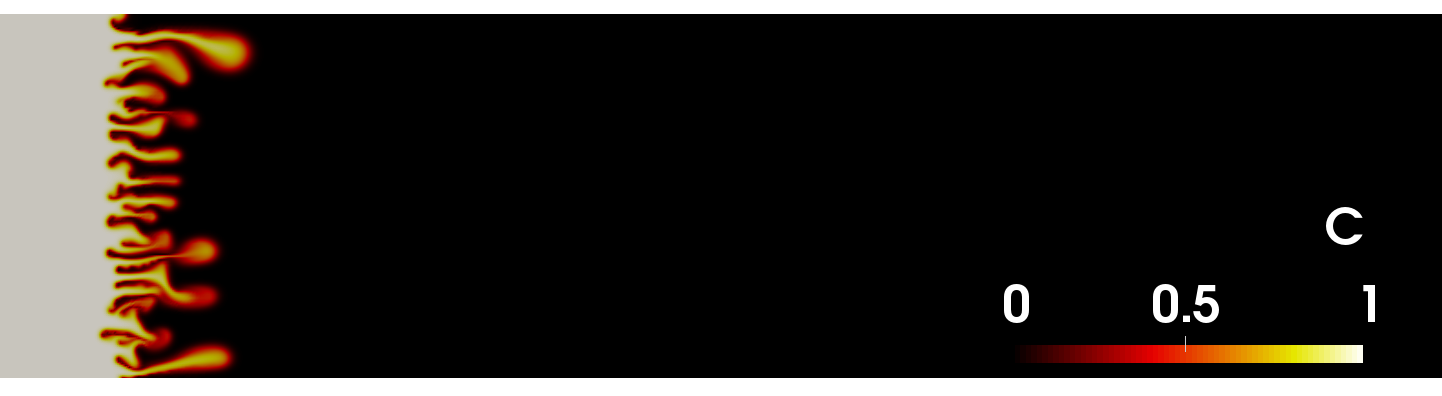}}\\
\subfloat[t=3]
{\includegraphics[width=0.45\textwidth]{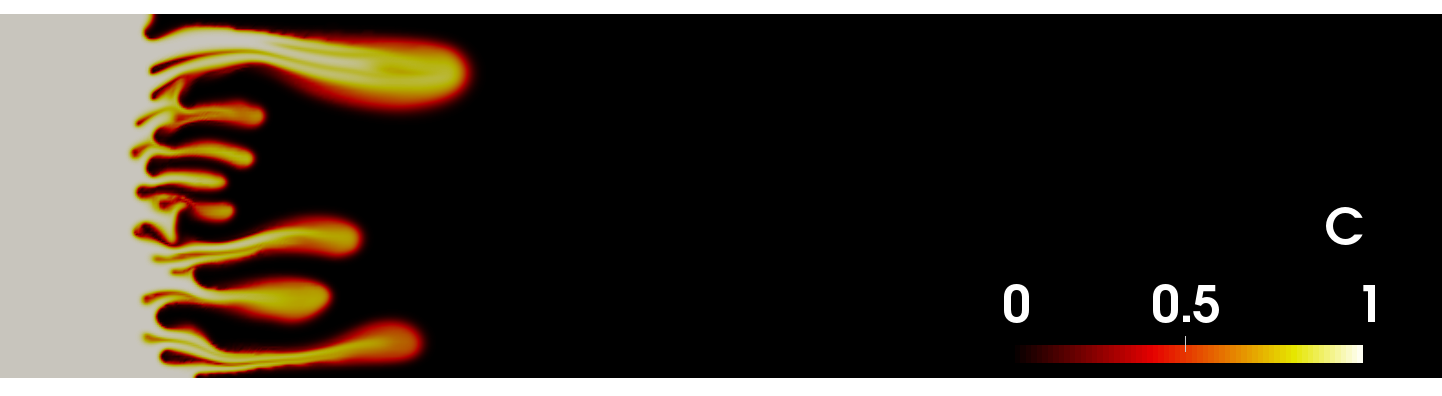}}
\subfloat[t=4]
{\includegraphics[width=0.45\textwidth]{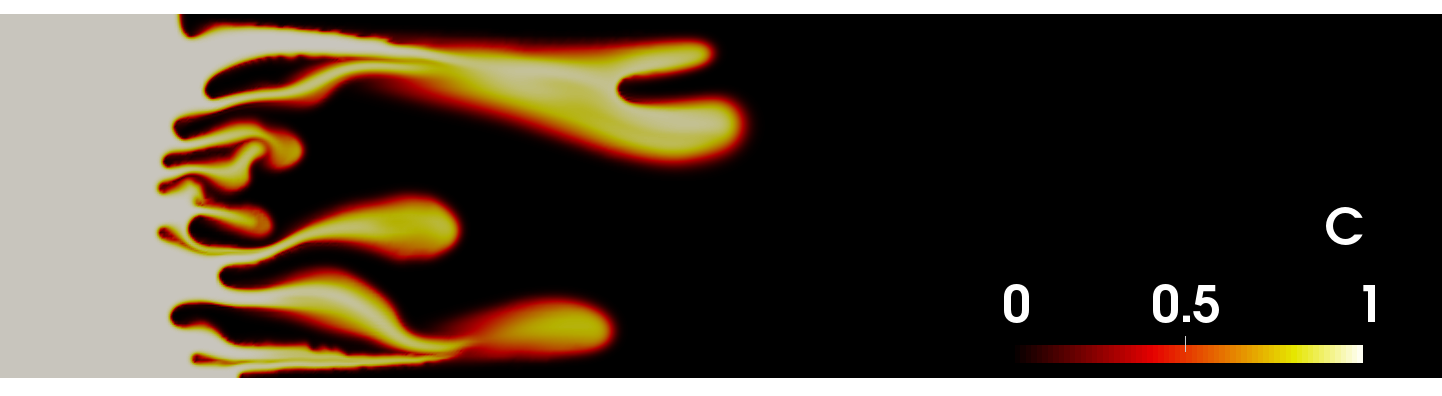}}
\caption{Example 5. Case 3 ($M=50$). Illustrates the concentration values for each time $t$.  The existing fluid is replaced by the another fluid and we observe viscous fingering instabilities. Fingers merge and also split.  }
\label{fig:ratio_50}
\end{figure}
\begin{figure}[!h]
\centering
\subfloat[n=100]
{\includegraphics[width=0.45\textwidth]{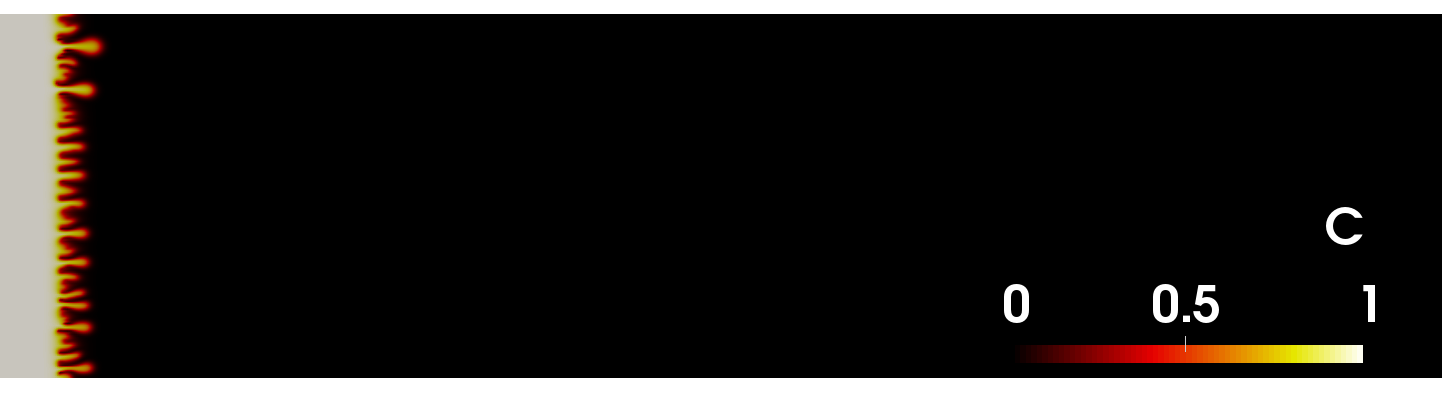}}
\hspace*{0.01in}
\subfloat[n=100]
{\includegraphics[width=0.45\textwidth]{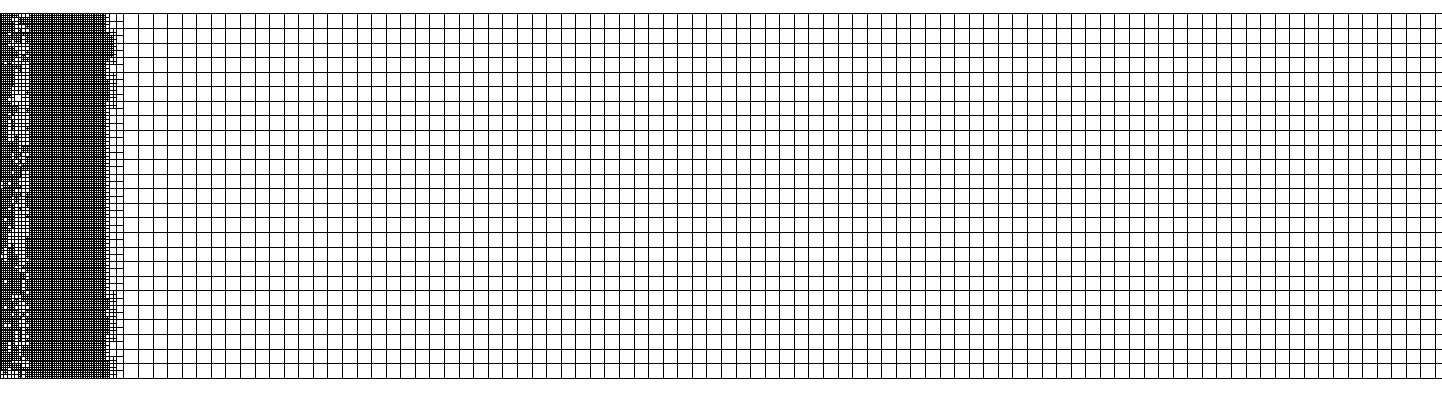}} \\
\subfloat[n=200]
{\includegraphics[width=0.45\textwidth]{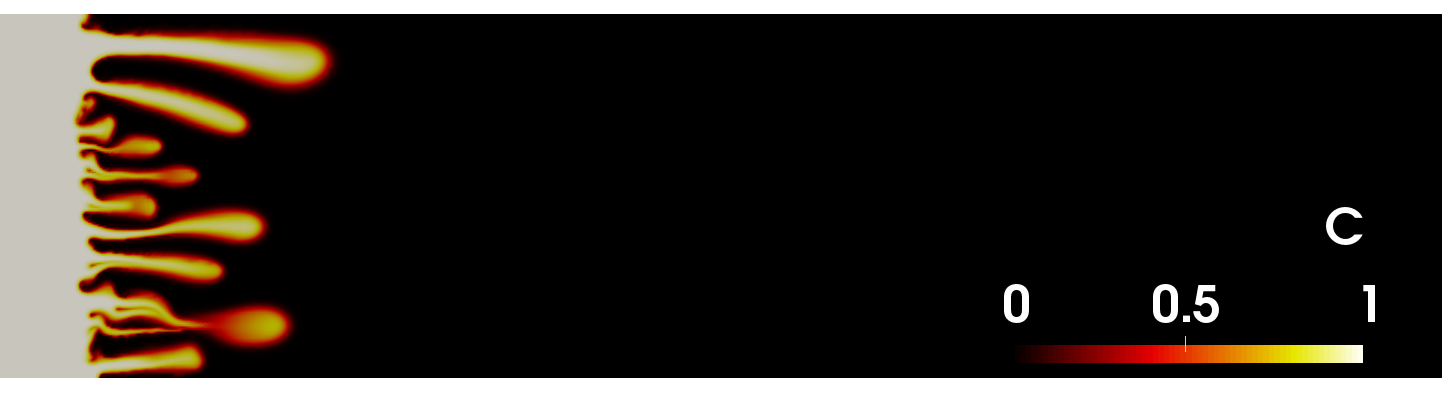}} 
\hspace*{0.01in}
\subfloat[n=200]
{\includegraphics[width=0.45\textwidth]{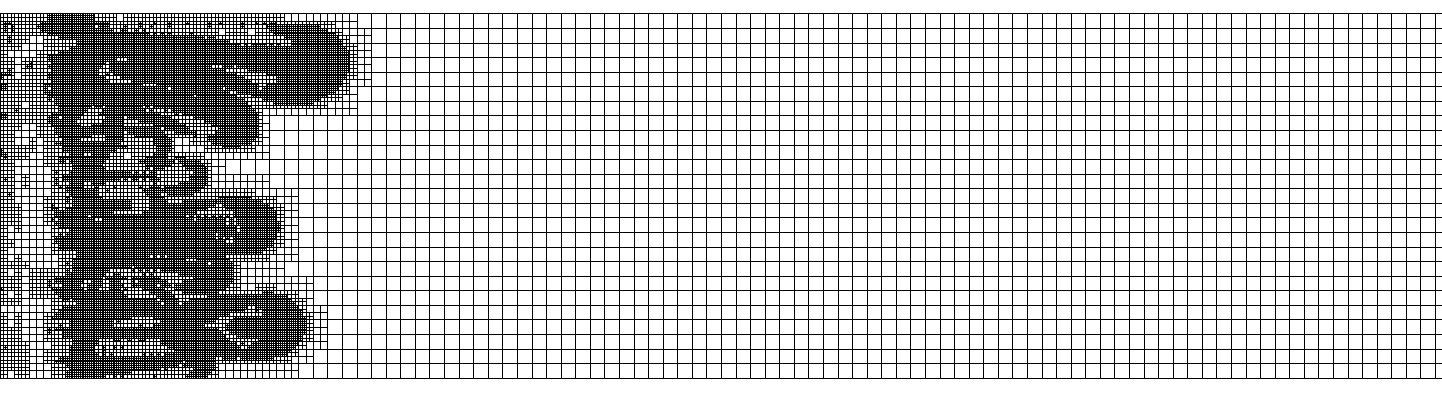}} \\
\subfloat[n=350]
{\includegraphics[width=0.45\textwidth]{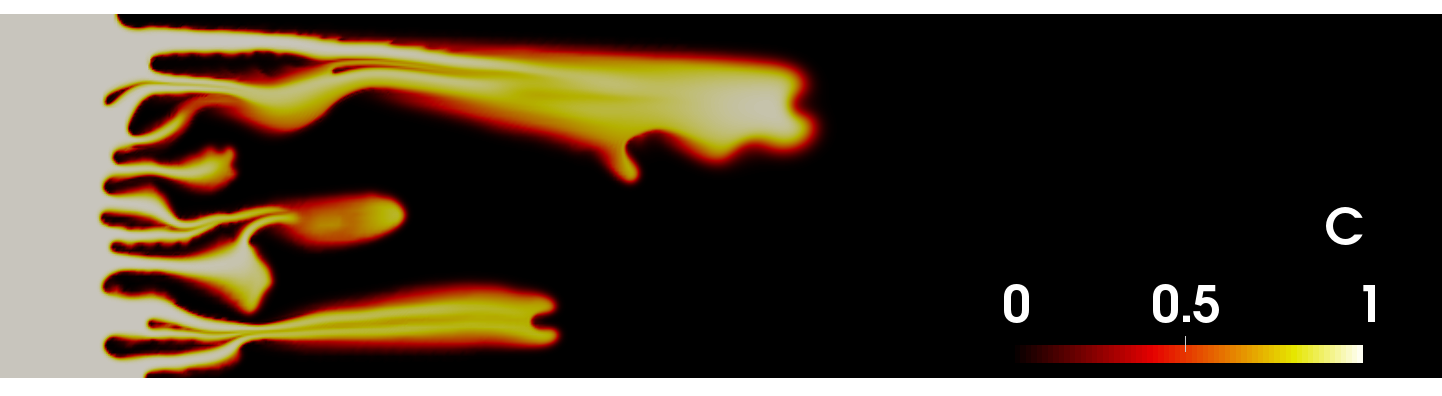}}
\hspace*{0.01in}
\subfloat[n=350]
{\includegraphics[width=0.45\textwidth]{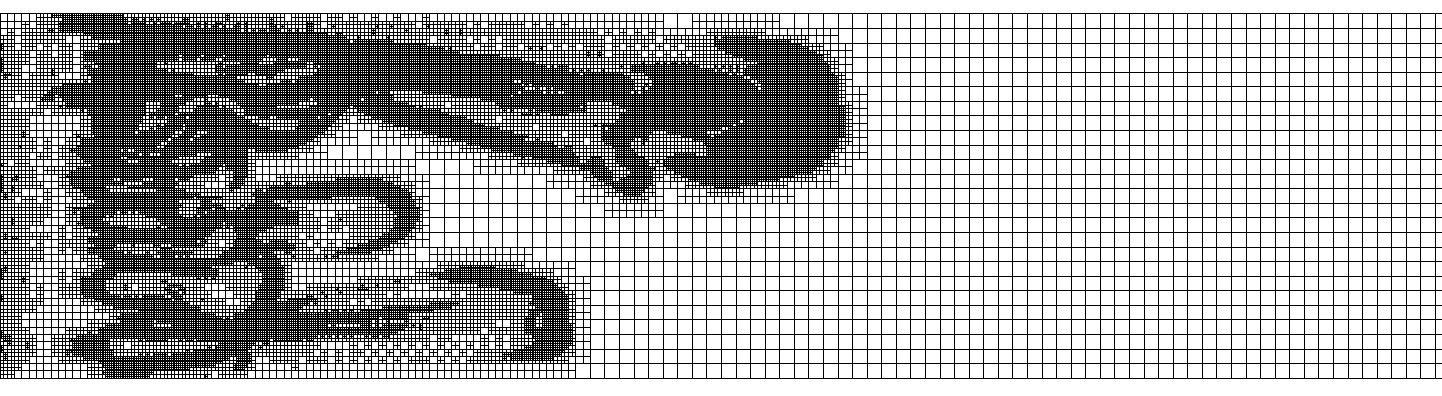}} \\
\subfloat[n=450]
{\includegraphics[width=0.45\textwidth]{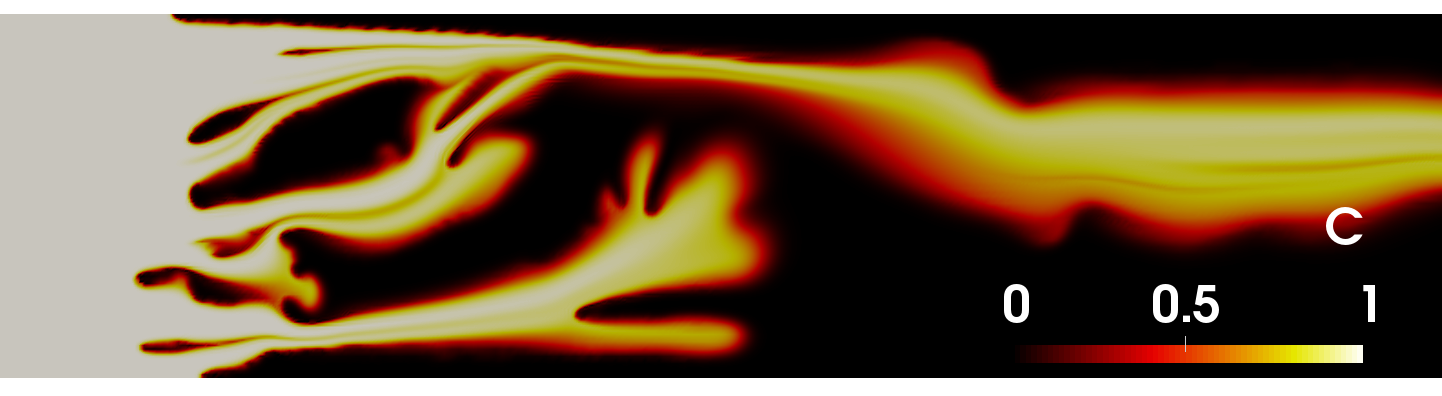}}
\hspace*{0.01in}
\subfloat[n=450]
{\includegraphics[width=0.45\textwidth]{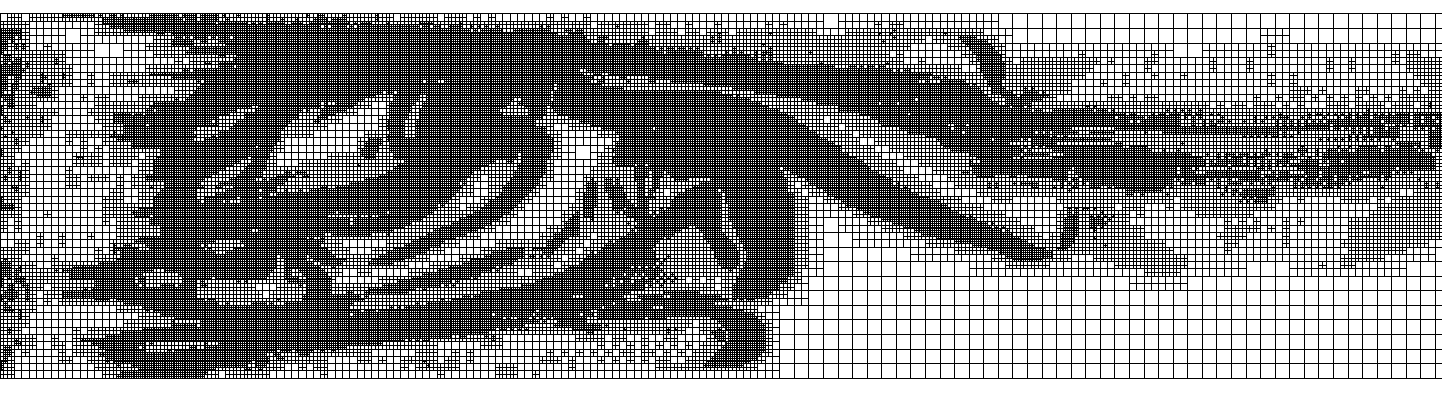}} 
\caption{Example 5. Case 6 ($M=150$). Concentration values for each time steps with corresponding adaptive mesh refinement at the right column. }
\label{fig:ratio_150}
\end{figure}

\begin{figure}[!h]
\centering
\subfloat[n=70]
{\includegraphics[width=0.45\textwidth]{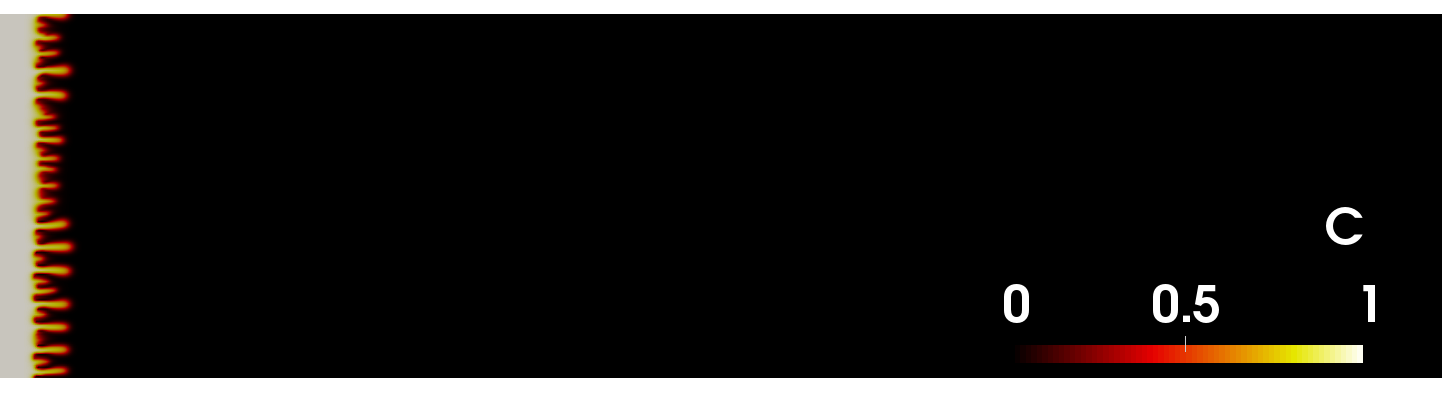}}
\hspace*{0.01in}
\subfloat[n=100]
{\includegraphics[width=0.45\textwidth]{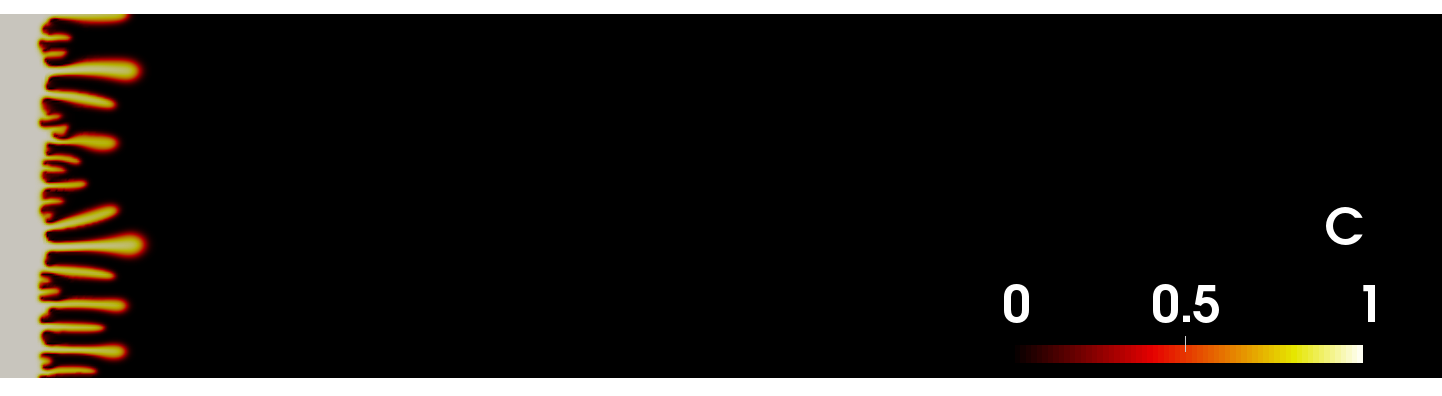}} \\
\subfloat[n=200]
{\includegraphics[width=0.45\textwidth]{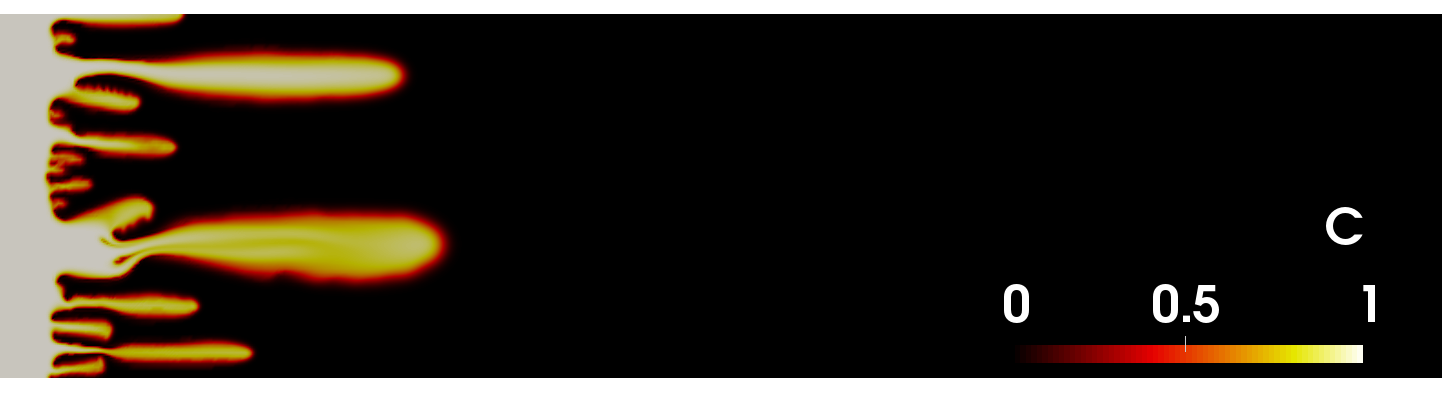}}
\hspace*{0.01in}
\subfloat[n=300]
{\includegraphics[width=0.45\textwidth]{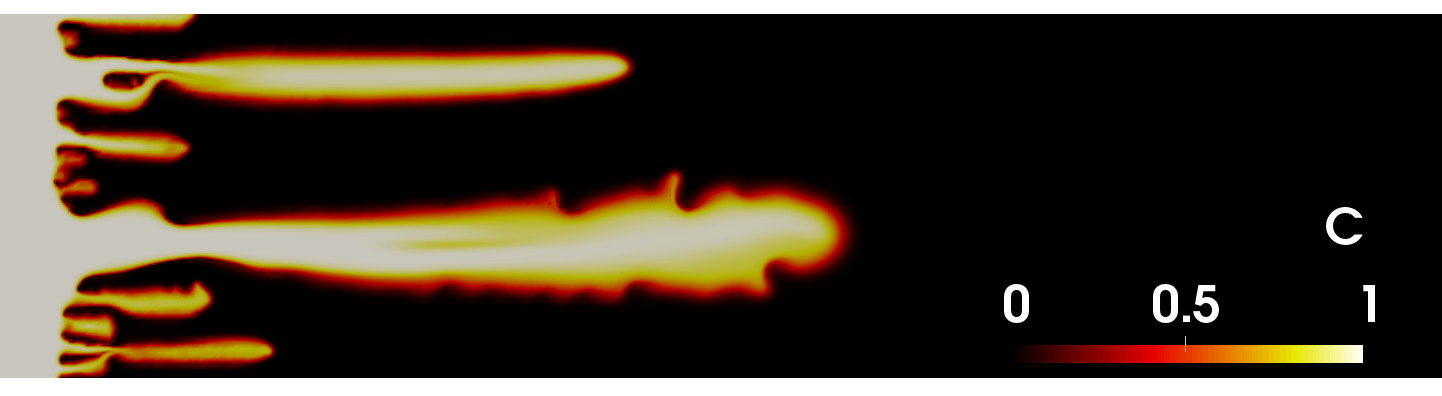}}
\caption{Example 5. Case 10 ($M=750$). Concentration values for each time step. We observe that the higher viscosity ratio could create more unstable fingers at more early time than lower viscosity ratios. See Figure \ref{fig:fingers}. }
\label{fig:ratio_750}
\end{figure}
\begin{figure}[!h]
\centering
\subfloat[Contour values of $C=0.5$]
{\includegraphics[scale=0.35]{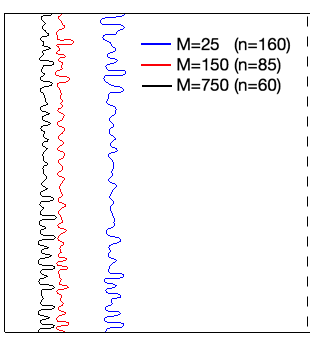}}
\subfloat[Time when the fingers are initiated versus viscosity ratio.  ]
{\includegraphics[scale=0.35]{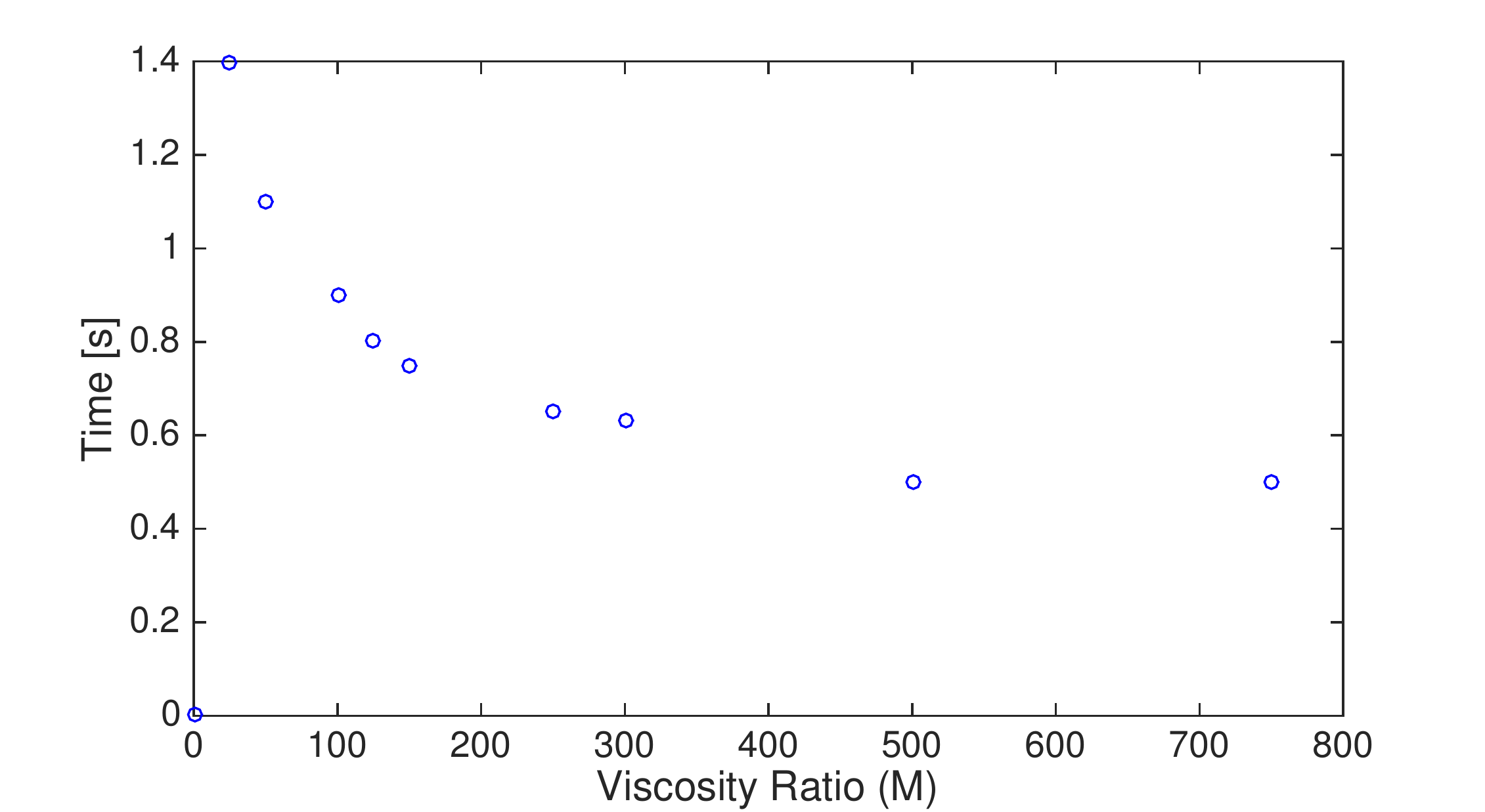}}
\caption{(a) illustrates the comparison of the fingers for each different viscosity ratios ($M$). We observe that the higher ratio initiates fingers at earlier time compare to the lower ratios. (b) shows the time when the fingers are initiated for each different viscosity ratios.  }
\label{fig:fingers}
\end{figure}

First, Figure \ref{fig:ratio_1} illustrates the Case 1 ($M=1$), where $\mu_s = \mu_0  = \SI{0.001}{\pascal\second}$ for both fluids. 
We do not observe any instability or fingering but the previous fluid is smoothly displaced by the incoming fluid. 
However, we do observe viscous fingering for all the remaining cases in Table 1, where $M>1$. 
We selectively illustrate three different cases,  $M = 50, 150$, and $M=750$, at Figures \ref{fig:ratio_50}, \ref{fig:ratio_150}, and \ref{fig:ratio_750}, respectively. 
Especially, in Figure \ref{fig:ratio_150} with $M=150$  we emphasize the effects of dynamic mesh refinement on capturing  viscous fingering instabilities. 

Is it well known that the behavior and the number of fingers varies with viscosity ratio, P\'{e}clet number, aspect ratio of the domain, and dispersion values \cite{moissis1993simulation}.
In early time, we observe a large number of small fingers followed by some of these initial fingers growing faster than others (e.g see Figure \ref{fig:ratio_50} (b) and (c)). 
At later time, the smaller fingers tend to merge with the larger ones as well as some of the fingers splitting (see Figure \ref{fig:ratio_50} (d)).
In addition, Figure \ref{fig:fingers} (a)-(b) illustrate that the higher viscosity initiates fingers at an earlier time than the lower viscosity ratios.
In addition, we observe that the fingers grow with a similar speed at lower viscosity ratios but there are some fingers that grow extremely faster than others  at high viscosity ratios (see Figure \ref{fig:ratio_750} (d)). 
{ These phenomena were also studied in \cite{saffman1958penetration}}.

\begin{figure}[!h]
\centering
\includegraphics[scale=0.35]{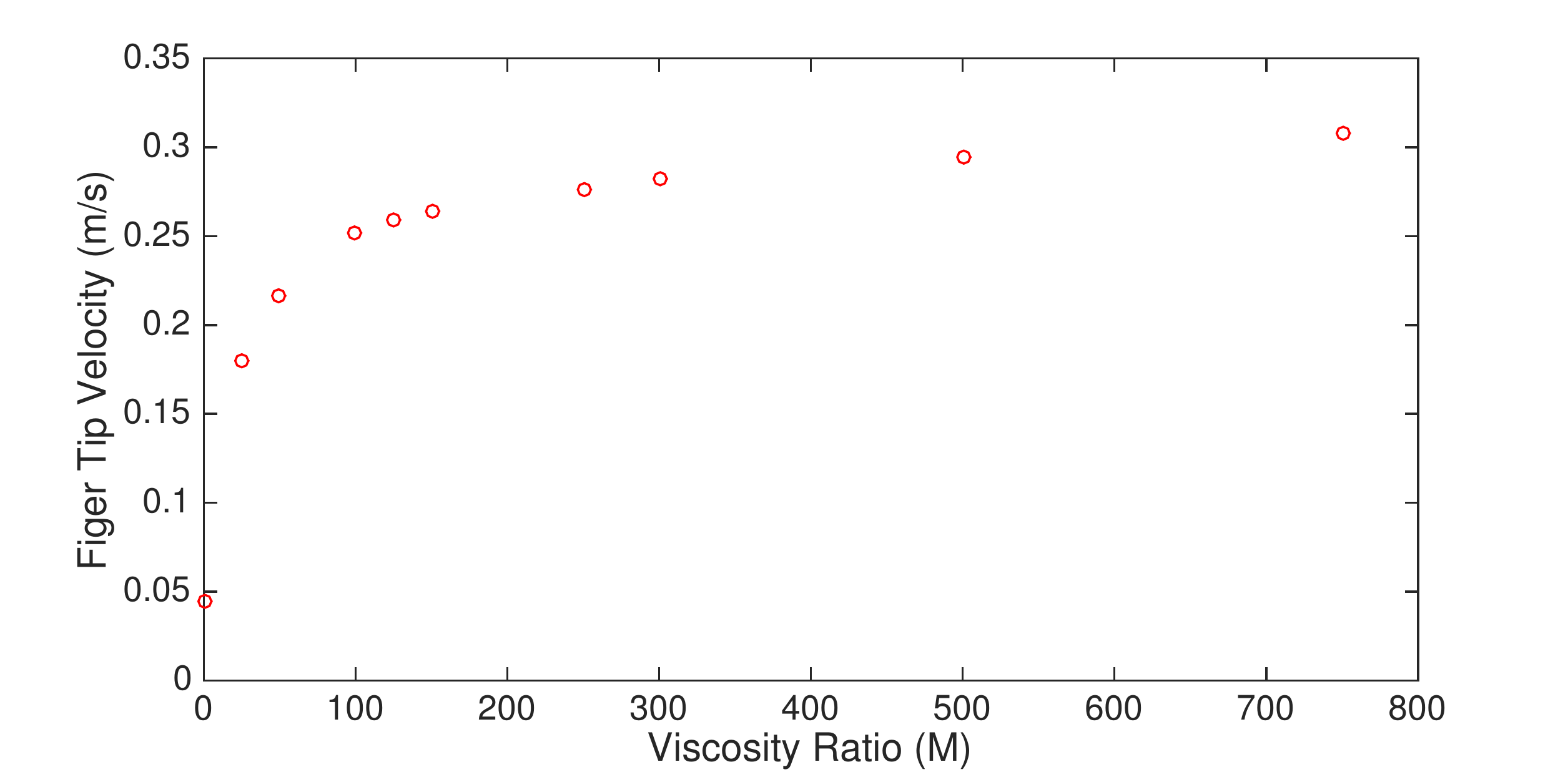}
\caption{Finger tip velocity versus the viscosity ratio.
The speed increases rapidly for smaller $M<300$, but it becomes almost plateau later. This is similar result as shown in the physical experiment \cite{malhotra2015experimental}.
}
\label{fig:finger_speed}
\end{figure}

Previously, 
petroleum engineers were interested in preventing
viscous fingering to improve oil production.
However, recently, viscous fingering is being employed for fracture propagation \cite{malhotra2014proppant, Lee2016} to efficiently transport the proppant to the tip of fracture.
Thus it has become important to study finger tip velocity \cite{malhotra2015experimental}.
{ 
Growth rate of the instabilities based on the linearization theory \cite{Homsy:1987cm,tan1986stability,mikelic91,CosBen,LaMaRaSaYo,BalRas}  has been studied in \cite{wooding1969growth,MenonOtto2005,manickam1995fingering,MenonOtto2006,YortsosSalin2006}. 
For example, \cite{YortsosSalin2006,MenonOtto2005} demonstrate that the leading finger tip velocity is bounded by $(M -1)^2/(M \ln M)$, where M is the end-point viscosity ratio.}
Figure \ref{fig:finger_speed} illustrates the finger tip velocity for each different viscosity ratios. 
The speed of the finger tip increases rapidly for smaller $M$, but it becomes almost a plateau for larger $M$. This behavior is very similar as observed in the recent experiment in (\cite{malhotra2015experimental}, Figure 14).

\subsubsection{Three dimensional Hele-Shaw cell}
In this section, we consider the three dimensional domain $\Omega = (0,0,0) \times (\SI{1}{\metre}, \SI{0.25}{\metre}, \SI{0.05}{\metre})$ to show the computational capabilities of our algorithm.  The numerical parameters are set to $h_{\min}=0.01$ and $\delta t = 0.01$ and other conditions are the same as in the previous section.
The maximum number of degrees of freedom for EG transport is approximately $325,000$ with the maximum cell number of $150,000$. 
Here we illustrate the case 4 with the viscosity ratio $M=100$ and the inflow is given by $p_{\text{in}}=0.01$. Figure \ref{fig:fingers_3D} illustrates the viscous fingering in a three dimensional Hele-Shaw cell with rectilinear flow.
\begin{figure}[!h]
\centering
\subfloat[n=100]
{\includegraphics[width=0.43\textwidth]{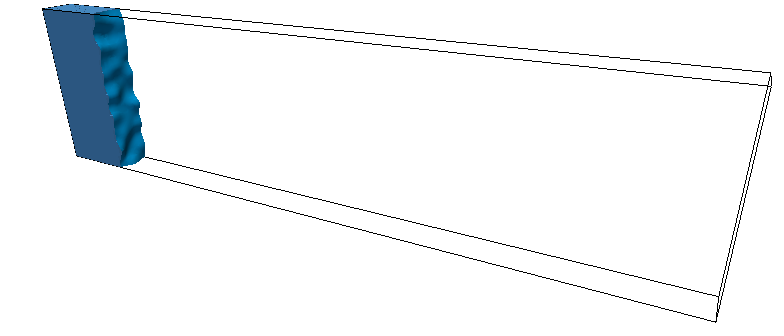}}
\hspace*{0.01in}
\subfloat[n=150]
{\includegraphics[width=0.43\textwidth]{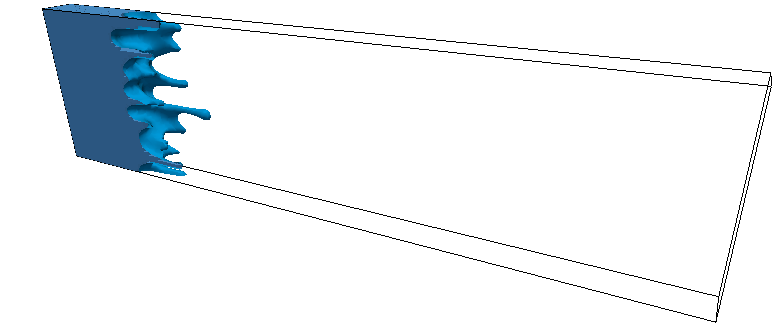}} \\
\subfloat[n=200]
{\includegraphics[width=0.43\textwidth]{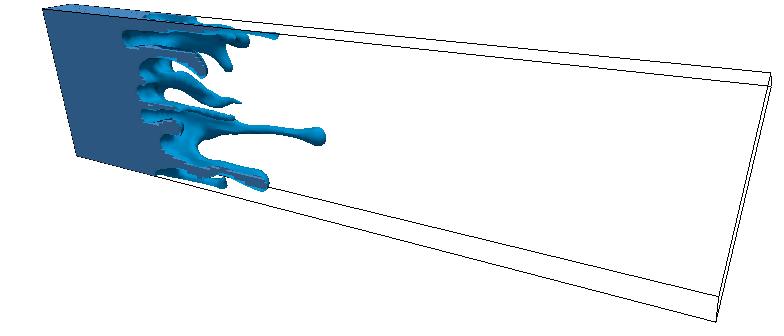}}
\hspace*{0.01in}
\subfloat[n=230]
{\includegraphics[width=0.43\textwidth]{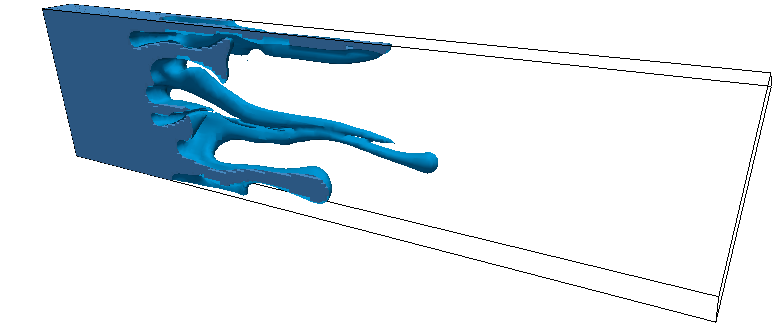}} \\
\subfloat[n=260]
{\includegraphics[width=0.43\textwidth]{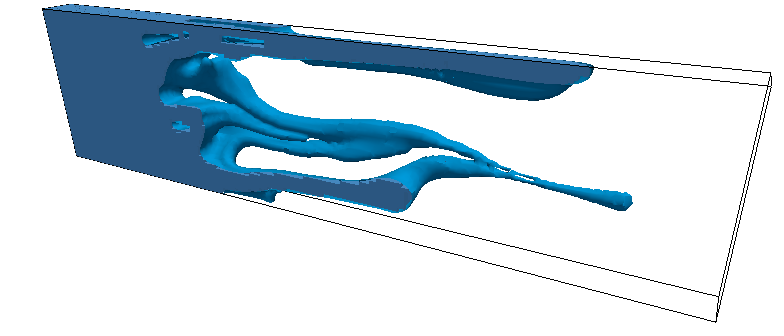}}
\hspace*{0.01in}
\subfloat[n=270]
{\includegraphics[width=0.43\textwidth]{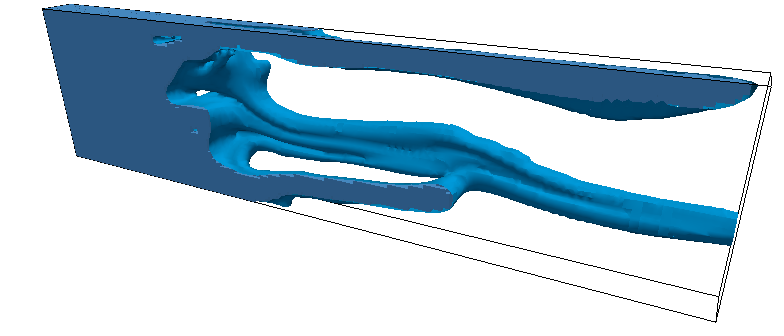}} 
\caption{Example 5. Case 4 ($M=100$) in three dimensional domain. Illustrates concentration values ($c <0.5$) for each time step.}
\label{fig:fingers_3D}
\end{figure}

\subsection{Hele-Shaw cell: Viscous fingering with a radial source flow.}
\label{sec:ex:radial}

In this final example, we present a radial flow displacement problem in a Hele-Shaw cell \cite{Homsy:1987cm,Homsy1987}, which was recently studied experimentally in \cite{bischofberger2014fingering}.
The computational domain is given as $\Omega = (\SI{0}{\metre},\SI{1}{\metre})^2$, and the fluid is injected at the source point $(0.5,0.5)$ with given $q / \rho_0=100$ and $c_q= 1$. 
The initial conditions for the pressure and concentration are set to zero, i.e
\begin{equation}
c_0 = 0  \  \text{ and  }  \  p_0  =0.
\end{equation}
The boundary conditions for the pressure and transport system are set to 
\begin{equation}
\bu \cdot \bn = 0 \quad 
\text{ and } \quad 
(-\bD(\bu)\nabla c) \cdot \bn = 0 
\text{ on } \partial \Omega \times (0,\mathbb{T}],
\end{equation}
respectively. In addition, the physical coefficient for diffusion and dispersion tensor are chosen as 
\begin{equation}
d_m = 1.8e^{-8} m^2/s,
\alpha_l = 1.8e^{-5} m^2/s, \mbox{ and }
\alpha_t = 1.8e^{-6} m^2/s.
\end{equation}
The numerical parameters for this example are set to $\delta t= 0.005$ and $h_{\min} = 0.01$ with stability coefficients 
$\lambda_{\textsf{Lin}} = \lambda_{\textsf{Ent}} = 1$. 
For adaptive mesh refinement, we set $\textsf{R}_{\max} = 9$ and  $\textsf{R}_{\min} = 7$.

Figure \ref{fig:radial_finger} illustrates the viscous fingering in Hele-Shaw cell by radial injection.  
Here $\mu_s = \SI{1}{\pascal \second}$ and $\mu_0 = \SI{0.001}{\pascal \second}$ with the viscosity ratio $M=1000$.
The finger tips split at later time as presented in  physical experiments
\cite{Homsy:1987cm,bischofberger2014fingering,Homsy1987}. 

\begin{figure}[!h]
\centering
\subfloat[$n=50$]
{\includegraphics[width=0.225\textwidth]{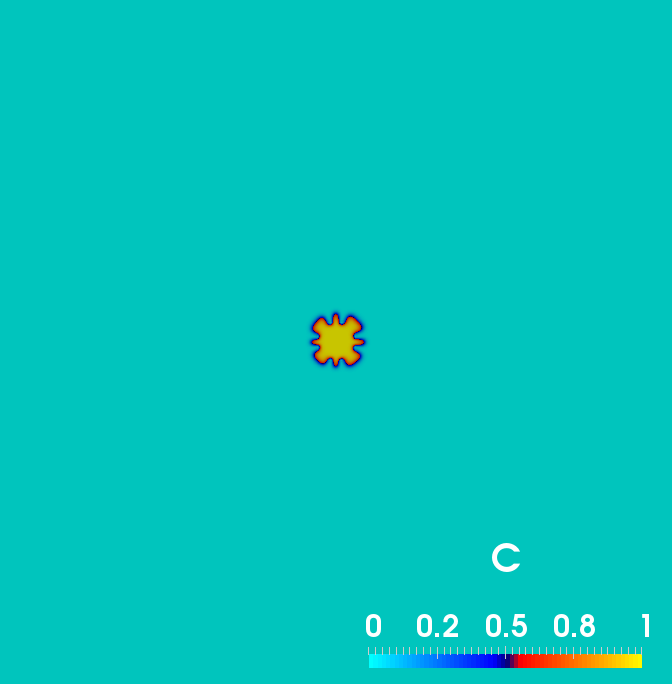}} 
\hspace*{0.05in}
\subfloat[$n=500$]
{\includegraphics[width=0.225\textwidth]{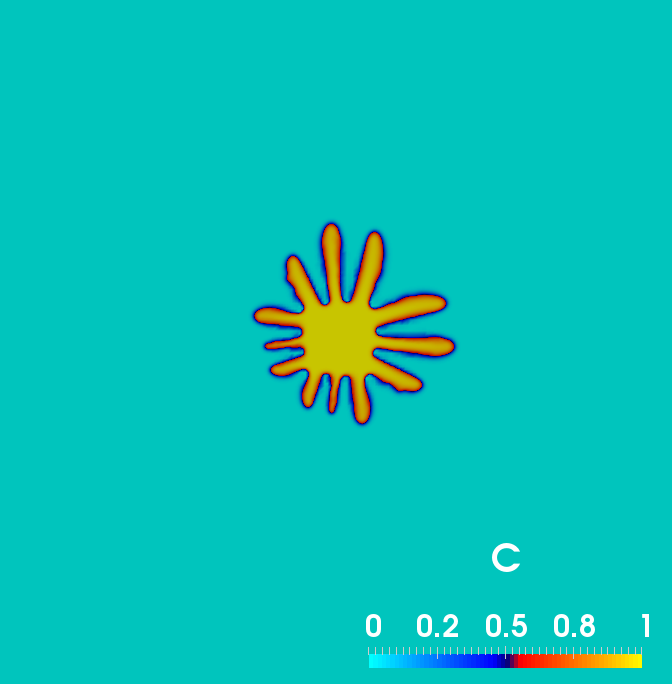}} 
\hspace*{0.05in}
\subfloat[$n=1000$]
{\includegraphics[width=0.225\textwidth]{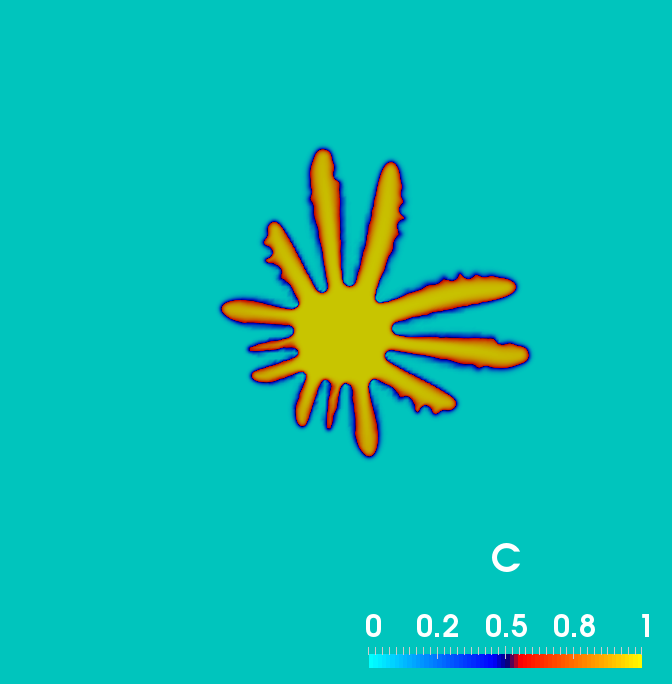}} 
\hspace*{0.05in}
\subfloat[$n=1880$]
{\includegraphics[width=0.225\textwidth]{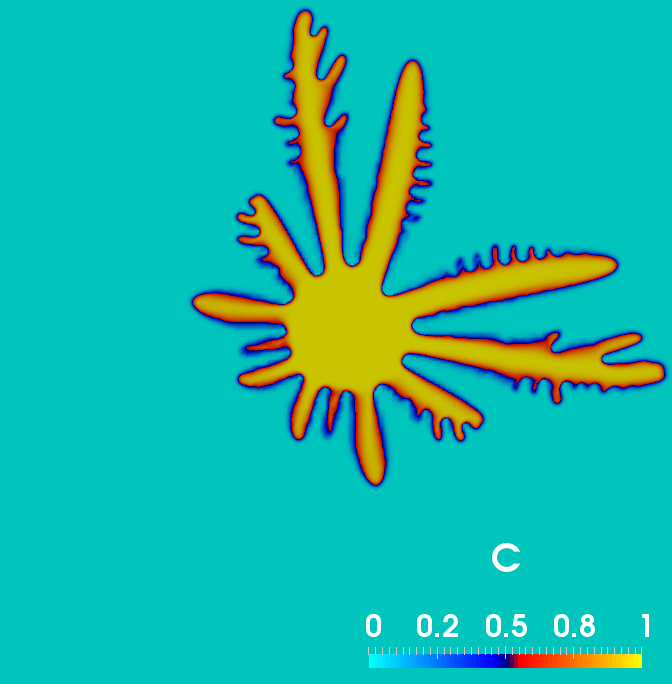}} 
\caption{Evolution of the fingers with viscosity ratio $M=1000$. The finger tip splits at later time as shown in some physical experiments \cite{Homsy:1987cm,bischofberger2014fingering,Homsy1987}. }
\label{fig:radial_finger}
\end{figure}

\section{Conclusion}
\label{sec:conclusion}
In this paper, we presented enriched Galerkin (EG) approximations for miscible displacement problems in porous media and Hele-Shaw cells. EG preserves local and global conservation for fluxes but has much fewer degrees of freedom compare to that of DG. For a higher order EG transport system, entropy residual stabilization is applied to avoid any spurious oscillations. 
In addition, dynamic mesh adaptivity employing entropy residual as an error indicator saves computational cost for large-scale computations.
Several examples including viscous fingering 
were constructed in order to demonstrate the performance of the algorithm.
This work can be extended for general flow (e.g Stokes)   including two phase flow.

\section*{Acknowledgments}
The research by S. Lee and  M. F. Wheeler was partially supported by 
{a} DOE grant DE-FG02-04ER25617,
{a} Statoil grant STNO-4502931834, and 
{an} Aramco grant UTA 11-000320.

\bibliographystyle{elsarticle-num} 
\bibliography{xyz}

\end{document}